\documentclass[preprint]{elsarticle}

\usepackage{graphicx}
\usepackage{xfrac,amsmath,amsmath}
\usepackage{esint,amsfonts,amsmath,amssymb,epsfig,mathrsfs,bm,accents}
\usepackage{ifthen,epic,eepic,color,stmaryrd,yfonts,mathtools}

\journal{Journal of \LaTeX\ Templates}

\numberwithin{equation}{section}
 \newtheorem{thm}{Theorem}[section]
 \newtheorem{lem}[thm]{Lemma}
 \newdefinition{rmk}[thm]{Remark}
 \newtheorem{definition}[thm]{Definition}
 \newtheorem{proposition}[thm]{Proposition}
 \newtheorem{corollary}[thm]{Corollary}
 \newproof{pf}{Proof}
 \newproof{pot}{Proof of Theorem \ref{thm2}}

\newcommand{\R}{\mathbb{R}}
\newcommand{\Rn}{{\R}^n}

\newcommand{\eps}{{\varepsilon}}

\newcommand{\diam}{{\text {diam}}}
\newcommand{\dist}{{\textup {dist}}}

\renewcommand\d{{\textup{d}}\,}

\newcommand\supp{{\rm supp}\,}
\newcommand\res{\mathop{\hbox{\vrule height 7pt width .5pt depth 0pt
\vrule height .5pt width 6pt depth 0pt}}\nolimits}

\newcommand{\cG}{{\mathcal{G}}}

\newcommand{\cH}{{\mathcal{H}}}

\newcommand{\cM}{{\mathcal{M}}}

\newcommand{\cC}{{\mathcal{C}}}
\newcommand{\cT}{{\mathcal{T}}}

\newcommand{\Pe}{{\mathscr{P}}}

\newcommand\F{{\mathbb F}}
\newcommand\N{{\mathbb N}}

\def\I#1{{\mathcal{A}}_{#1}}

\newcommand{\Iq}{{\mathcal{A}}_Q}
\def\a#1{\left\llbracket{#1}\right\rrbracket}

\newcommand{\de}{\partial}

\newcommand{\etaa}{{\bm{\eta}}}

\newcommand{\graph}{\text{graph}}
\newcommand{\loc}{\textup{loc}}

\newcommand\Sing{\textup{Sing}}

\renewcommand{\d}{\textup{d}}
\renewcommand{\Sing}{\textup{Sing}}

\newcommand\sV{{\mathscr V}}
\newcommand\sC{{\mathscr C}}
\newcommand\sW{{\mathscr W}}
\newcommand\sZ{{\mathscr Z}}

\newenvironment{itemizeb}
{\begin{itemize}\itemsep=2pt}{\end{itemize}}

\begin{document}

\begin{frontmatter}

\title{Improved estimate of the singular set of Dir-minimizing $Q$-valued functions
via an abstract regularity result}

\author{Matteo Focardi\corref{cor1}}
\ead{focardi@math.unifi.it}
\address{DiMaI, Universit\`a degli Studi di Firenze, Viale Morgagni~67/A, 50134 Firenze, Italy}

\author{Andrea Marchese}
\ead{marchese@mis.mpg.de}

\author{Emanuele Spadaro}
\ead{spadaro@mis.mpg.de}
\address{Max-Planck-Institut f\"ur Mathematik in den Naturwissenschaften, Inselstrasse~22, 04103 Leipzig, Germany}

%
\cortext[cor1]{Corresponding author}
%

\begin{abstract}
In this note we prove an abstract version of a recent quantitative stratification priciple introduced by 
Cheeger and Naber (\textit{Invent. Math.}, 191 (2013), no. 2, 321--339; \textit{Comm. Pure Appl. Math.}, 
66 (2013), no. 6, 965--990).
Using this general regularity result paired with an $\eps$-regularity theorem we provide a new estimate of 
the Minkowski dimension of the set of higher multiplicity points of a Dir-minimizing $Q$-valued function.
The abstract priciple is applicable to several other problems: we recover recent results in the literature and
we obtain also some improvements in more classical contexts.
\end{abstract}


\end{frontmatter}


\section{Introduction}

\paragraph{{\bf An abstract regularity result}}

We propose an abstraction of a quantitative stratification principle introduced and developed in a series of papers 
by Cheeger and Naber \cite{ChNa1, ChNa2}, Cheeger, Haslhofer and Naber \cite{ChHaNa1, ChHaNa2} 
and Cheeger, Naber and Valtorta \cite{ChNaVal}.

The interest in finding general formulations of this kind of regularity results is driven by
a number of important applications in geometric analysis.
Apart from those contained in the papers quoted above, we mention the cases of Dir-minimizing
$Q$-valued maps according to Almgren, of varifold with bounded mean curvature and of almost minimizers of the perimeter. 
The former is treated in details in  \S~\ref{s:Q} and \S~\ref{s:dirminstrat}, the latters in \S~\ref{s:vrfld}.
We explicitly remark that the papers \cite{ChHaNa1, ChHaNa2} deal also with parabolic examples, a case that is not covered by our results.

To our knowledge the first example in this direction of abstraction is the general regularity theorem 
proven by Simon \cite[Appendix A]{Sim83} based on the so called \textit{dimension reduction argument} introduced 
by Federer in his pioneering work \cite{Fed2}. Similarly, the paper by White \cite{Wh} generalizes the 
refinement of Federer's reduction argument made by Almgren in his big regularity paper \cite{Alm}.

The basic principle and the main ingredients of our abstract formulation
can be explained roughly as follows.
\smallskip

\begin{quotation}
\noindent\textsc{Abstract stratification:} 
\noindent\textit{the set of points where a solution to a geometric problem is faraway at every scale from being 
homogeneous with $k+1$ indipendent invariant directions has Minkowski  dimension less than or equal to $k$.}
\end{quotation}

\smallskip

The main sets of quantities we consider are:
\begin{itemizeb}
\item[(a)] a family of density functions $\Theta_s$, increasing w.r.t.~$s\geq 0$; 

\item[(b)] a family of distance functions $\d_k$, $k\in\{0, \ldots, m\}$, measuring 
the distance from $k$-invariant homogeneous solutions. 

\end{itemizeb}

In addition, we assume suitable compatibility conditions, namely

\begin{itemizeb}

\item[(i)] a quantitative differentiation principle that allows to 
quantify the number of those scales for which closeness to homogeneous solutions fails, 
and that tipically follows in the applications from monotonicity type formulas;

\item[(ii)] a consistency relation between the distances $\d_k$: if a solution is close to a $k$-invariant one 
and additionally is $0$-invariant with respect to another point away from the invariant $k$-dimensional space, 
then it is actually close to a  $(k+1)$-invariant solution (see \S~\ref{sh} for the detailed formulation).
\end{itemizeb}

This set of hypotheses is common to many problems in geometric analysis such as the
Dirichlet minimizing multiple valued 
functions, harmonic maps, almost minimizing currents and several others (see \cite{ChNa1}-\cite{ChNaVal} for other applications).
Indeed, the stratification result and the estimate on the Minkowski dimension only 
depend on assumptions (i) and (ii), thus making the common aspects of all previous results clear.

It  turns out that there is a simple connection between White's approach to Almgren's stratification and the one outlined above. 
In \S~\ref{s:bianco} we show how to recast the result by White in our framework. 
In this respect, we stress that the stratification in \cite{ChNa2}
and in our Theorem~\ref{t:nostra stratiFICA} can be applied to some cases
not covered by the ideas in \cite{Wh}, such as stationary harmonic maps
(cp.~\cite[Corollary~2.6]{ChNa2}, \S~\ref{ss:shm} and \cite[Section~6]{Wh}).

\smallskip

Our main application of the abstract stratification principle is outlined in the following paragraph.

\paragraph{{\bf Application to $Q$-valued functions}}

In the regularity theory for higher codimension minimal surfaces 
(in the sense of \textit{mass minimizing integer rectifiable currents})
a fundamental role is played by the multiple valued functions introduced by Almgren
in \cite{Alm}, which turn out to be the correct blowup limits for the
analysis of singularities (see also \cite{DSmemo,DSannali, DS1, DS2, DS3}
for a simplified new proof of the result in \cite{Alm}).

Following \cite{DSmemo}, a $Q$-valued function $u$ is a measurable map
from a bounded open subset $\Omega\subset \R^n$ (for simplicity
we always assume that the boundary of $\Omega$ is smooth) taking values
in the space of positive atomic measures in $\R^m$ with mass $Q$, namely
\[
\Omega\ni x \mapsto u(x) \in \Iq(\R^m) := \left\{ \sum_{i=1}^Q \a{p_i} \; : \; p_i \in \R^m \right\},
\]
where $\a{p}$ denotes the Dirac delta at $p$.
Almgren proves in \cite{Alm} (cp.~also \cite{DS3}) that the blowups
of higher codimension mass minimizing integral currents are actually
graphs of $Q$-valued functions $u$ in a suitable Sobolev class
$W^{1,2}(\Omega, \Iq(\R^m))$ minimizing a generalized Dirichlet energy
(cp.~\cite[Definition~0.5]{DSmemo}):
\[
\int_{\Omega} |Du|^2 \leq \int_{\Omega} |Dv|^2 \quad \forall\; v \in W^{1,2}(\Omega, \Iq(\R^m)),\;v\vert_{\de\Omega} = u\vert_{\de\Omega},
\]
(explicit examples of Dir-minimizing $Q$-valued functions are given in \cite{Spadaro}).

In order to estimate the size of the singular set of a minimizing current
it is essential to bound the dimension of the set of points
where the graph of a Dir-minimizing $Q$-valued function has higher multiplicity.
Almgren's main result in the analysis of multiple valued functions is in fact
an estimate of the Hausdorff dimension of the set 
$\Delta_Q$ of multiplicity $Q$ points of a Dirichlet
minimizing $Q$-valued function $u$, i.e.~the set of points $x \in \Omega$ such that 
$u(x) = Q\a{p}$ for some $p \in \R^m$,
which turns out not to exceed $n-2$ in the case it does not coincide with $\Omega$
(cp. \cite[Proposition~3.22]{DSmemo}).

In this paper we improve Almgren's result by showing an estimate of
the Minkowski dimension of $\Delta_Q$. 
To this aim we denote by $\cT_{r}(E) := \{ z \in \R^n : \dist(z, E) < r\}$
the tubular neighborhood of radius $r$ of a given set $E\subset\R^n$.

\begin{thm}\label{t:A}
Let $u :\Omega \to \Iq(\R^m)$ be a Dir-minimizing function, where $\Omega\subset \R^n$
is a bounded open set with smooth boundary.
Then either $\Delta_Q = \Omega$, or for every $\Omega'\subset \subset \Omega$
the Minkowski dimension of $\Delta_Q \cap \Omega'$ is less than or equal to $n-2$, i.e.~for every $\Omega'\subset \subset \Omega$ and for every $\kappa_0 > 0$
there exists a constant $C>0$ such that
\begin{equation}\label{e:introQ}
|\cT_r(\Delta_Q\cap \Omega')| \leq C\, r^{2 - \kappa_0} \quad \forall \;0<r<\dist(\Omega', \de \Omega).
\end{equation}
\end{thm}

We also obtain a stratification result for the whole set of singular points of multiple valued functions that, even if known to the experts, we were not able to find in the literature.
To this aim we introduce the following notation. Given a $Q$-valued function
$u:\Omega \to \Iq(\R^m)$, we denote by $\Sing_u \subset \Omega$ its singular set, i.e.~$x_0\not\in \Sing_u$ if and only if there exists $r>0$ such that
\[
\graph(u\vert_{B_r(x_0)}) := \{(x, y)\in \R^{n \times m}\;:\; |x-x_0| < r,\, y \in \supp(u(x))\}
\]
is a smooth $n$-dimensional embedded submanifold (not necessarily connected).
For every $k\in\{0,\ldots, n\}$, we define the subset $\Sing^k_u$ of the singular set $\Sing_u$
made of those points having all tangent functions with at most $k$ independent directions of 
invariance (we refer to \S~\ref{ss:B} for the precise definition).

\begin{thm}\label{t:B}
Let $u :\Omega \to \Iq(\R^m)$ be a Dir-minimizing function, where $\Omega\subset \R^n$
is a bounded open set with smooth boundary, and let
$\Sing^k_u$ be the singular strata defined in \S~\ref{ss:B}.
Then, $\Sing_u = \Sing^{n-2}_u$ and
\begin{gather}
\Sing^0_u \quad \text{is countable} \label{e:B1}\\
\dim_\cH(\Sing^k_u) \leq k \qquad \forall\; k\in\{1, \ldots, n-2\}.\label{e:B2}
\end{gather}
\end{thm}

In the case $Q=2$ a more refined
analysis by Krummel and Wickramasekera \cite{KrWi} shows the rectifiability
of the singular set, remarkably improving Almgren's work.

We prove Theorems~\ref{t:A} and \ref{t:B} as a consequence of our abstract stratification principle.
More precisely, Theorem~\ref{t:B} is a direct consequence of it, while Theorem~\ref{t:A} requires a further 
stability property deduced by an $\eps$-regularity result (see Proposition~\ref{e.eps reg minchione}).

\paragraph{{\bf Applications to generalized submanifolds}}

In the final section \S~\ref{s:vrfld} we apply the abstract stratification principle to varifolds with bounded mean 
curvature and almost minimizers of the perimeter, two relevant cases for applications that are not covered by the results 
in \cite{ChNa2}. Also in these cases we derive some improvements of well-known estimates for the singular set.
Stratification for the singular set of stationary varifolds with bounded mean curvature is addressed in 
\S~\ref{s:vrfld tub neg}. Eventually, in Theorem~\ref{t:am minchioschi} we give a bound on the Minkowski dimension 
of the singular set of an almost minimizer of the perimeter rather than the classical  Hausdorff dimension estimate, 
and in Theorem~\ref{t:higher int} we show higher integrability for its generalized second fundamental form.

\paragraph{{\bf On the organization of the paper}}
A few words are worthwhile concerning the structure of the paper. 
The first two sections of the paper are devoted to the abstract regularity results. 
In particular, \S~\ref{s:abstract} contains the estimate of the volume of the tubular neighborhood of the singular strata 
given in Theorem~\ref{t:tub-neig} (which is proved in the first part of \S~\ref{proofthm}) and the abstract stratification 
in Theorem~\ref{t:nostra stratiFICA}. In order to make 
our statements and hypotheses recognizable and ``natural'' to the readers, we illustrate them in \S~\ref{modelex} for the model examples of area minimizing currents and harmonic maps. 
The last part of \S~\ref{proofthm} is devoted to the comparison with the results by White in \cite{Wh}. 
Then, we specialize our results to the case of $Q$-valued functions in \S\ref{s:dirminstrat}, the needed preliminaries 
are collected in \S~\ref{s:Q}. We finally focus on varifolds with suitable hypotheses on their mean curvature and on 
almost minimizers of the perimeter in \S~\ref{s:vrfld}.

%
%
\section*{Acknowledgements}
For this research E.~Spadaro has been partially supported by GNAMPA
Gruppo Nazionale per l'Analisi Matematica, 
la Probabilit\`a e le loro Applicazioni of the Istituto Nazionale di Alta Matematica (INdAM) through 
a Visiting Professor Fellowship. 
E.~Spadaro is very grateful to the DiMaI ``U. Dini'' of the University of Firenze 
for the support during the visiting period.

Part of this work was conceived when M. Focardi was visiting the Max Planck Institut in Leipzig. 
He would like to warmly thank the Institute for providing a very stimulating scientific atmosphere
and for all the support received.

\section{Abstract Stratification}\label{s:abstract}

The general abstract approach we propose is based on two main sets of quantities:
namely, a family of density functions $\Theta_s$ and 
an increasing family of distance functions
$\d_k$.

\subsection{Densities and distance functions}\label{pa:1.1}
Let $\Omega\subset\R^n$ be open and bounded,
and for every $s\geq 0$ set $\Omega^s:= \{x \in \Omega : \dist(x,\de \Omega) \geq 2s\}$.
We assume the following.

\begin{itemize}
\item[(a)] For every $s$ such that $\Omega^s\neq\emptyset$, there exist functions
$\Theta_s\in L^\infty(\Omega^s)$
such that
\[
0\leq \Theta_s(x) \leq \Theta_{s'}(x),
\]
for all $0 \leq s<s'$ and for all $x \in \Omega^{s'}$.
Moreover, for every $s_0>0$ there exists $\Lambda_0=\Lambda_0(s_0)>0$ such that
\[
\Theta_s(x) \leq \Lambda_0,
\]
for every $0\leq s\leq s_0$ and for every $x\in\Omega^{s_0}$.

\item[(b)]  Setting $U := \{(x,s): x\in \Omega^s, \, \Theta_0(x) > 0 \}$,
there exist a positive integer $m\leq n$ and control functions
$\d_k:U \to [0,+\infty)$ for  $k\in\{0,\ldots, m\}$
such that
\[
\d_0\leq \d_1 \leq \cdots \leq \d_m.
\]
\end{itemize}

\subsection{Structural Hypotheses}\label{sh}
These two sets of quantities are then related by the following
structural hypotheses.
\begin{itemize}
\item[(i)] For every $s_0>0$, $\eps_1>0$
there exist $0<\lambda_1(s_0,\eps_1),\,\eta_1(s_0,\eps_1) <\sfrac{1}{4}$ such that
if $(x,s) \in U$, with $x\in\Omega^{s_0}$ and $s<s_0$, then
\[
\Theta_{s}(x)-\Theta_{\lambda_1s}(x)\leq\eta_1
\quad\Longrightarrow\quad
\d_0(x,s) \leq \eps_1.
\]
\item[(ii)] For every $s_0>0$, for every $\eps_2, \tau\in (0,1)$ there exists 
$0<\eta_2(s_0,\eps_2,\tau)\leq \eps_2$ such that if $(x,5s)\in U$, with $x\in\Omega^{s_0}$ 
and $5s<s_0$, satisfies for some $k \in \{0, \ldots, m-1\}$
\[
\d_k(x,4s) \leq \eta_2 \quad \text{and}\quad \d_{k+1}(x,4s) \geq \eps_2,
\]
then
there exists a $k$-dimensional linear subspace $V$ for which
\[
\d_0(y,4s)> \eta_2 \quad \forall\; y \in B_{s}(x)\setminus
\cT_{\tau s}(x+V),
\]
where $\cT_{\tau s}(x+V) := \{ z : \dist(z, x+V) < \tau s\}$
is the tubular neighborhood of $x+V$ of radius $\tau s$.
\end{itemize}

\subsection{Volume of the neighborhoods of singular strata}

The sets we consider in our estimates are the following.

\begin{definition}[Singular Strata]
For every $0<\delta<1$, $0<r\leq r_0$
and for every $k\in \{0, \ldots, m-1\}$ we set
\begin{equation}
\mathcal{S}^k_{r, r_0, \delta} :=
\big\{x\in \Omega^{r_0} :\Theta_0(x)>0\quad\text{and}\quad \d_{k+1}(x,s) \geq \delta \quad
\forall\; r\leq s \leq r_0 \big\}
\end{equation}
and
\begin{equation}
\mathcal{S}^k_{r_0, \delta}:=\bigcap_{0<r \leq r_0}
\mathcal{S}^k_{r,r_0, \delta}
\quad\text{and}\quad
\mathcal{S}^k_{r_0} := \bigcup_{0<\delta<1} \mathcal{S}^k_{r_0, \delta}.
\end{equation}
\end{definition}
Note that, by the monotonicity of the control functions,
$\mathcal{S}^k_{r, \delta}\subset \mathcal{S}^{k'}_{r', \delta'}$
if $\delta'\leq \delta$, $r\leq r'$ and $k\leq k'$.

\medskip

Our abstract stratification result relies on the following
estimate for the tubular neighborhoods of the singular strata. Its proof is postponed to \S \ref{proofthm}.

\begin{thm}\label{t:tub-neig}
Under the Structural Hypotheses in \S~\ref{sh},
for every $\kappa_0,\delta \in (0,1)$ and $r_0>0$
there exists $C=C(\kappa_0, \delta, r_0, n,\Omega)>0$ such that
\begin{equation}\label{e:tubulare}
\vert \cT_r(\mathcal{S}^k_{r, r_0, \delta})\vert\leq C\, r^{n-k-\kappa_0}
\quad \forall \; 0 < r < r_0
\quad\forall\;k\in\{1, \ldots, m-1\}
\end{equation}
\begin{equation}\label{e:countable}
 \mathcal{S}^0_{r_0, \delta} \text{ is countable}.
\end{equation}
\end{thm}

\subsection{Hausdorff dimension of the singular strata}
It is now an immediate consequence of Theorem~\ref{t:tub-neig} the following stratification result.

\begin{thm}\label{t:nostra stratiFICA}
Under the Structural Hypotheses in \S~\ref{sh} for every $r_0>0$ the estimate 
$\dim_\cH(S^{k}_{r_0}) \leq k$ holds for $k\in\{1, \ldots, m-1\}$. Moreover, $S^0_{r_0}$ is countable.
\end{thm}

\begin{pf}
Indeed Theorem~\ref{t:tub-neig} implies that $\dim_{\cM}(S^k_{r_0,\delta}) \leq k$, where $\dim_{\cM}$ is the Minkowski dimension.
Since the Hausdorff dimension of a set is always less than or equal to the Minkowski dimension, we also infer that
\[
\dim_{\cH}(S^k_{r_0}) \leq
\dim_{\cH}\Big(\bigcup_{\delta>0} S^{k}_{r_0,\delta}\Big) \leq k
\]
because, being the union monotone, it is enough to consider a countable set of parameters.
\end{pf}

\subsection{Minkowski dimension of the singular strata}\label{e:deltadep}
The dependence of the constant $C$ in \eqref{e:tubulare} on $\delta$ prevents the derivation of an estimate on the Minkowski dimension of the singular strata $S^k_{r_0}$.
Nevertheless, if such dependence drops, then Theorem~\ref{t:tub-neig} turns actually into an estimate on the Minkowski dimension of the singular strata 
which is not implied by the Almgren's stratification principle.

\begin{thm}\label{t:minchioschi}
Under the hypotheses of Theorem~\ref{t:tub-neig}, if for some $\delta_0>0$ and $k\in \{0, \ldots, m-1\}$
\begin{equation}\label{e:stability}
\mathcal{S}^k_{r_0,\delta} = \mathcal{S}^k_{r_0} \quad \forall \; \delta \in (0, \delta_0), 
\end{equation}
then for every $0<\kappa_0 <1$ and $r_0>0$
there exists $C=C(\kappa_0, \delta_0, r_0, n,\Omega)>0$ such that
\begin{equation}\label{e:tubular}
\vert \cT_r(\mathcal{S}^k_{r_0})\vert\leq C\, r^{n-k-\kappa_0}
\quad \forall \; 0 < r < r_0.
\end{equation}
In particular $\dim_\cM(\mathcal{S}^k_{r_0}) \leq k$.
\end{thm}

\subsection{Examples}\label{modelex}
The meaning of the Structural Hypotheses in \S~\ref{sh} is very well illustrated by 
the two familiar examples of area minimizing currents and stationary
harmonic maps treated in \cite{ChNa2} for which Theorem~\ref{t:tub-neig} and \ref{t:nostra stratiFICA} hold.
Moreover for area minimizing currents of codimension one in $\R^n$ Theorem~\ref{t:minchioschi} can be also applied for $k= n-8$.

\subsubsection{Area minimizing currents}\label{pa.mc}
Let $T$ be an $m$-dimensional area minimizing integral current in $\Omega$.
Then we can set
\begin{equation*}
\Theta_s(x):=\frac{\|T\|(B_s(x))}{w_m s^m} \quad \text{for }\;s>0
\quad \text{and} \quad \Theta_0(x) := \lim_{r\downarrow 0^+} \Theta_r(x)
\end{equation*}
and for $k\in\{0, \ldots, m\}$
$$
\d_k(x,s):=\inf\big\{\mathbb{F}\big((T_{x,s}-C)\res B_1\big): C\text{ is $k$-conical \& area minimizing} \big\},
$$
where
\begin{itemizeb}
\item $T_{x,r}$ is the rescaling of the current around any
point $x \in \R^n$ at scale $r>0$: 
\begin{equation}\label{e.rescalings}
T_{x,r} := \big(\eta_{x,r}\big)_\# T
\end{equation}
and the push-forward is done via the proper map $\eta_{x,r}$ given by $y\mapsto \sfrac{(y-x)}{r}$;
\item  $\mathbb{F}$ is the flat norm (see \cite[\S~31]{Sim83});
\item an $m$-dimensional current $C$ in $\R^n$ is \emph{$k$-conical}
for $k \in \{0, \ldots, m\}$, if there 
exists a linear subspace $V \subset \R^n$ of dimension bigger than or equal to $k$ such that 
\[
\text{$T_{x,r} = T$ for all $r>0$ and $x\in V$.}
\]
Note that a $0$-conical current is simply a cone with respect to the origin.
\end{itemizeb}
One can choose $\Lambda_0(r_0):=\sfrac{\mathbb{M}(T)}{\omega_m r_0^m}$. Then (a) is a consequence of the Monotonicity Formula (see \cite[Theorem~17.6]{Sim83}) and (b) follows 
from the inclusion of $k$-conical currents in the $k'$-conical ones when $k'\leq k$.
With this choice, the structural hypoteses in \S~\ref{sh} are satisfied, indeed (i) is an other consequence of the Monotonicity Formula and (ii) follows from a rigidity property of cones sometimes called ``cylindrical blowup'' 
(see \cite[Lemma 35.5]{Sim83}).

Then the quantitative stratification principle in Theorem~\ref{t:tub-neig}
recovers the corresponding result in \cite{ChNa2}:

\smallskip

\begin{quotation}
\noindent\textit{the set of points that are faraway
from $(k+1)$-conical area minimizing currents,
at every scale in $[r, r_0]$,
has Minkowski dimension less than or equal to $k$.}
\end{quotation}

\subsubsection{Stationary harmonic maps}\label{ss:shm}
Similarly let $u\in W^{1,2}(\Omega,\mathscr{N})$ be a stationary harmonic
map from an open set $\Omega\subset\R^n$ to a Riemannian manifold $(\mathscr{N}^m,h)$ 
isometrically embedded in some Euclidean space $\R^p$ (see, e.g., \cite{Sim-verde}).
We can set 
\[
\Theta_s(x):=s^{2-n}\int_{B_s(x)}|\nabla u|^2dy,\quad s\in\big(0,\dist(x,\partial\Omega)\big),
\]
and for every $k\in\{0,\dots,n\}$ \[
\d_k(x,r):=\inf_{v\in\sC_k}
\fint_{B_1}\dist^2_{\mathscr{N}}\big(u_{x,r},v\big)dy,
\]
with 

\begin{itemizeb}
\item $u_{x,r}(y):=u(x+ry)$ for $x\in\Omega$ and $r\in\big(0,\dist(x,\partial\Omega)\big)$;
\item a measurable map $v$ is said to be \emph{$k$-conical} if there exists a vector space $V$ with $\dim V\geq k$ that leaves $v$ invariant, i.e.
\begin{equation}\label{e:harmconic1}
v(x)=v(y+x)\quad\forall\,x\in\R^n,\,y\in V,
\end{equation}
and such that $v$ is $0$-homogeneous with respect to the points in $V$, i.e. 
\begin{equation}\label{e:harmconic2}
v(y+x)=v(y+\lambda\,x)\quad\forall\,x\in\R^n,\,y\in V\text{ and }\lambda>0;
\end{equation}
\item $\sC_k:=\{\text{$v:B_1\to\mathscr{N}$ $k$-conical}\}$ .
\end{itemizeb}

Assumption (a) in \S~\ref{pa:1.1} is easily verified and the monotonicity
formula
\begin{equation*}
\Theta_r(x)-\Theta_s(x)=\int_s^r\int_{\partial B_t(x)}t^{2-n}\left|\frac{\partial u}{\partial t}\right|^2d\cH^{n-1}dt
\end{equation*}
together with an elementary contradiction argument
show that the Structural Hypothesis (i) in \S~\ref{sh} is satisfied.
Moreover the Structural Hypothesis (ii) follows similarly to the one
for minimizing currents (cp.~\cite{ChNa2} for more details), thus
leading to the stratification of Theorem~\ref{t:tub-neig}.

\medskip

In Section~\ref{s:vrfld} we give other applications of this abstract
regularity result to the case of varifolds with bounded variation
and almost minimizers of the mass in codimension one.

%
%

\section{Proof of the Abstract Stratification and comparison with Almgren's Stratification}\label{proofthm}
To begin with, we state a simple consequence of the Structural Hypothesis (ii) (cp. \S~\ref{sh}) 
in the following
\begin{lem}\label{l:cattura}
For every $s_0>0$, for every $\eps, \tau\in (0,1)$ there exists $0<\gamma_0\leq \eps$
such that if $(x,5s)\in U$, with $x\in\Omega^{s_0}$ and $5s<s_0$, satisfies for some $k \in \{0, \ldots, m-1\}$
\[
\d_0(x,4s) \leq \gamma_0 \quad\text{and}\quad \d_{k+1}(x,4s) \geq \eps,
\]
then there exists
a linear subspace $V$ with $\dim(V)\leq k$ such that
\begin{equation}\label{e:inclusione}
y \in B_{s}(x) \; \textup{ \& }\; \d_0(y,4s) \leq \gamma_0
\quad\Longrightarrow\quad
y \in \cT_{\tau s}(x+V).
\end{equation}

\end{lem}

\begin{pf}
Let $\gamma_0 \leq \gamma_1 \leq \ldots\leq \gamma_{k+1}$ be set as
$\gamma_{k+1}=\eps$ and $\gamma_{j-1} = \eta_2(s_0,\gamma_j, \tau)$ with $\eta_2$ the
constant in the Structural Hypothesis (ii).
Let $i\in \{0,\ldots, k\}$ be the smallest index such that
$\d_{i+1}(x, 4s) \geq \gamma_{i+1}$ (which exists because of
the assumption $\d_{k+1}(x,4s) \geq \eps=\gamma_{k+1}$).
Then, applying the Structural Hypothesis (ii) we deduce that
there exists an $i$-dimensional linear subspace $V$ such that every point $y \in B_s(x)$
with $\d_0(y,4s) \leq \gamma_0\leq \gamma_{i}$ belongs to the tubular neighborhood
$\cT_{\tau s}(x+V)$.
\end{pf}

In the proof of Theorem~\ref{t:tub-neig} we shall repeatedly use the following elementary covering argument.

\begin{lem}\label{l.ricoprimento banana}
For every measurable set $E \subset \R^n$ with finite measure and for every $\rho>0$, there exists a finite covering $\{B_\rho(x_i)\}_{i \in I}$ of $\cT_{\sfrac{\rho}{5}}(E)$ with $x_i \in E$ and
\begin{equation}\label{e.stima banana}
\cH^0(I) \leq \frac{5^n\,|\cT_{\sfrac{\rho}{5}}(E)|}{\omega_n\, \rho^n}.
\end{equation}
\end{lem}

\begin{pf}
Consider the family of balls $\{B_{\sfrac{\rho}{5}}(x)\}_{x \in E}$. By the Vitali $5r$-covering lemma, there exists a finite subfamily $\{B_{\sfrac{\rho}{5}}(x_i)\}_{i \in I}$ of disjoint balls such 
that $\cT_{\sfrac{\rho}{5}}(E) \subset \cup_{i \in I}B_{\rho}(x_i)$. By a simple volume comparison we conclude \eqref{e.stima banana}.
\end{pf}

We are now ready to prove Theorem~\ref{t:tub-neig}.
\begin{pf}[of Theorem~\ref{t:tub-neig}]
We start fixing a parameter $\tau = \tau(n, \kappa_0)>0$ such that
\begin{equation}\label{e:tau}
\omega_n\,\tau^{\frac{\kappa_0}{2}} \leq 20^{-n}.
\end{equation}
We choose the other constants involved in the Structural Hypotheses in the following way:
\begin{enumerate}
\item let $\gamma_0 \leq \gamma_1 \leq \ldots\leq \gamma_k$ be such that
$\gamma_k=\delta$ and $\gamma_{j-1} = \eta_2(r_0, \gamma_j, \tau)$ for every $j\in\{1,\ldots,k\}$
as in the Structural Hypothesis (ii);
\item let $\lambda_1=\lambda_1(r_0,\gamma_0)$ and $\eta_1=\eta_1(r_0,\gamma_0)$ be
as in the Structural Hypothesis (i);
\item fix $q\in \N$ such that $\tau^q \leq \lambda_1$.
\end{enumerate}

We divide the proof into four steps.

\medskip

\noindent\emph{Step 1: reduction to dyadic radii.}
Let $\Lambda_0=\Lambda_0(r_0)$ given in \S~\ref{pa:1.1}. It suffices to prove \eqref{e:tubulare} for every $r$ of the form
$r=\frac{r_0\tau^p}{5}$ with $p\in \N$ such that $p \geq p_0 := q+M+1$ and $M:=\lfloor q\,\Lambda_0/\eta_1 \rfloor$.
Indeed for $\frac{r_0\tau^{p_0}}{5} < s < r_0$ we simply have
\begin{align*}
\vert \cT_s(\mathcal{S}^k_{s, r_0, \delta})\vert &\leq |\Omega| \leq 
\frac{|\Omega|}{\left(\frac{r_0\tau^{p_0}}{5}\right)^{n-k-\kappa_0}}\, s^{n-k-\kappa_0}\\
& = C_2(\kappa_0, \delta, r_0, n,\Omega) \, s^{n-k-\kappa_0}.
\end{align*}
On the other hand, if we assume that \eqref{e:tubulare} holds with a constant $C_1>0$ for every $r$ of the form $r=\frac{r_0\tau^p}{5}$ with $p \geq p_0$, we conclude
that for $r\tau < s < r$ it holds
\begin{align*}
\vert \cT_{s}(\mathcal{S}^k_{s, r_0, \delta})\vert & \leq \vert \cT_{r}(\mathcal{S}^k_{r, r_0, \delta})\vert \leq C_1\, r^{n-k-\kappa_0} \leq C_1\, \tau^{k+\kappa_0 - n}\, s^{n-k-\kappa_0}.
\end{align*}
Therefore setting $C := \max\{\tau^{k+\kappa_0 - n}\,C_1, \, C_2\}$ we deduce that \eqref{e:tubulare} holds for every $r \in (0, r_0)$.

\medskip

\noindent\emph{Step 2: selection of good scales.}
Fix a value $p\in \N$ with $p \geq p_0$ as above and set $r=\sfrac{r_0\tau^p}{5}$.
For all $(x,r_0)\in U$ we have
\begin{align*}
\sum_{l=q}^{p}
\Theta_{4\tau^{l}\,r_0}(x)-\Theta_{4\tau^{l+q}\,r_0}(x) &=
\sum_{l=q}^{p}\sum_{i=l}^{l+q-1}
\Theta_{4\tau^{i}\,r_0}(x)-\Theta_{4\tau^{i+1}\,r_0}(x) \\
&\leq q\sum_{h=q}^{p+q-1}\big(\Theta_{4\tau^{h}\,r_0}(x)-\Theta_{4\tau^{h+1}\,r_0}(x)\big)\\
&=q\,\big(\Theta_{4\tau^q r_0}(x)-\Theta_{4\tau^{p+q}\,r_0}(x) \big)
\leq q\,\Lambda_0.
\end{align*}
Therefore, there exist at most $M$ indices $l\in \{q,\ldots, p\}$ for which it does not hold that
\begin{equation}\label{e:boni in du cess}
\Theta_{4\tau^{l}\,r_0}(x)-\Theta_{4\tau^{l+q}\,r_0}(x) \leq \eta_1.
\end{equation}
For any subset $A \subset \{q,\ldots, p\}$ with cardinality
$M$ we consider
\[
S_{A} := \left\{x \in S^k_{r,r_0, \delta} : \eqref{e:boni in du cess} \;\text{ holds }\;\forall\; l \not\in A \right\}.
\]
We prove in the next step that
\begin{equation}\label{e:stima basica}
\vert \cT_r(S_{A})\vert\leq C\, r^{n-k-\frac{\kappa_0}2}
\end{equation}
for some $C=C(\kappa_0,\delta, r_0, n,\Omega)>0$.
From \eqref{e:stima basica} one concludes because
the number of subsets $A$ as above is estimated by
\[
\binom{p-q+1}{M}\leq (p-q+1)^{M} \leq C\,|\log r|^M
\]
for some $C(\kappa_0,\delta, r_0, n)>0$, and
\[
|\cT_r(S^k_{r,r_0, \delta})| \leq \sum_A |\cT_r(S_A)| \leq
C\,  |\log r|^M\, r^{n-k-\frac{\kappa_0}2}\leq
C\,r^{n-k-\kappa_0}
\]
for some $C(\kappa_0,\delta, r_0, n,\Omega)>0$.

\medskip

\noindent\emph{Step 3: proof of \eqref{e:stima basica}.}
We estimate the volume of $\cT_r(S_A)$ by covering it iteratively 
with families of balls centered in $S_A$ and with radii $\tau^jr_0$ for $j\in \{q, \ldots, p\}$.
We can then proceed as follows.
Firstly we consider a cover of $\cT_{\sfrac{\tau^q r_0}{5}}(S_A)$ made of balls
$\{B_{\tau^q r_0}(x_i)\}_{i\in I_q}$ with $x_i \in S_A$ and by a straightforward use of Lemma~\ref{l.ricoprimento banana}
\[
\cH^0(I_q) \leq 5^n\tau^{-nq}r_0^{-n}\big(\diam(\Omega)+1\big)^n.
\]
Iteratively, for every $j\in\{q+1,\ldots, p\}$, we assume to be given
the cover $\{B_{\tau^{j-1}r_0}(x_i)\}_{i\in I_{j-1}}$ of $\cT_{\sfrac{\tau^{j-1}r_0}{5}}(S_A)$,
and we select a new cover of $\cT_{\sfrac{\tau^{j}r_0}{5}}(S_A)$
which is made of balls of radii $\tau^{j}r_0$ centered in $S_A$
according to the following two cases:
\begin{itemizeb}
\item[(a)] $j-1 \in A$,
\item[(b)] $j-1 \notin A$.
\end{itemizeb}

\textit{Case \textup{(a)}.}
For every $x_i$ in the family at level $j-1$, using Lemma~\ref{l.ricoprimento banana}
we cover $S_A \cap B_{\tau^{j-1}r_0}(x_i)$
with finitely many balls $B_{\sfrac{\tau^{j}r_0}{2}}(y_l)$ with $y_l \in
S_A \cap B_{\tau^{j-1}r_0}(x_i)$ and the cardinality of the cover is bounded by
\[
\frac{5^n\,|B_{(\tau^{j-1} + \sfrac{\tau^{j}}{10})\,r_0}(x_i)|}{\omega_n \left(\sfrac{\tau^j r_0}{2} \right)^n} \leq 
20^n\,\tau^{-n}
\]
(note that $\cT_{\sfrac{\tau^jr_0}{10}}(S_A \cap B_{\tau^{j-1}r_0}(x_i)) \subset B_{(\tau^{j-1} + \sfrac{\tau^{j}}{10})\,r_0}(x_i)$).
We claim next that the union of $B_{\tau^{j}r_0}(y_l)$
covers the tubular neighborhood
\[
\cT_{\frac{\tau^jr_0}{5}}(S_A \cap B_{\tau^{j-1}r_0}(x_i)).
\]
Indeed for every $z \in \cT_{\sfrac{\tau^jr_0}{5}}(S_A \cap B_{\tau^{j-1}r_0}(x_i))$
there exists $z' \in S_A \cap B_{\tau^{j-1}r_0}(x_i)$ such that $|z-z'| < \sfrac{\tau^jr_0}{5}$.
Since $z' \in B_{\sfrac{\tau^{j}r_0}{2}}(y_l)$ for some $y_l$,
then $z \in B_{\tau^{j}r_0}(y_l)$.

Therefore, collecting all such balls, the cardinality of the new covering is estimated by
\begin{equation}\label{e:cardinality}
\cH^0(I_{j}) \leq 20^n\,\tau^{-n}\, \cH^0(I_{j-1}).
\end{equation}

\medskip

\textit{Case \textup{(b)}.}
If $j-1 \notin A$, then \eqref{e:boni in du cess} holds with $l=j-1$.
By the Structural Hypothesis (i) and the choice of $\lambda_1, \eta_1$ in (2) and $\tau$ in (3) at the beginning of the proof, we have that
$\d_0(x,4\tau^{j-1}r_0)\leq \gamma_0$ for every $x\in S_A$.
Since $x_i \in S_A \subset S^k_{r,r_0,\delta}$ we have also
$\d_{k+1}(x_i,4\tau^{j-1}r_0) \geq \delta$.
We can then apply Lemma~\ref{l:cattura} and conclude that
\begin{equation}\label{e:intrappolamento}
S_{A} \cap B_{\tau^{j-1}r_0}(x_i)
\subset \cT_{\tau^jr_0}(x_i+V)
\end{equation}
for some linear subspace $V$ of dimension less than or equal to $k$.
Note that
\begin{equation}\label{e:volume}
|\cT_{\tau^{j}r_0}((x_i+V)\cap B_{\tau^{j-1}r_0}(x_i))|
\leq \omega_n \, \tau^{n-k}\, |B_{\tau^{j-1}r_0}(x_i)|.
\end{equation}
Thus applying Lemma~\ref{l.ricoprimento banana} we find a covering of
$\cT_{\sfrac{\tau^jr_0}{5}}(S_{A})$ with balls $B_{r_0\tau^{j}}(y_l)$ such that $y_l \in S_A$ and using \eqref{e:volume} the cardinality of the covering is bounded by 
\begin{equation}\label{e:cardinality2}
\cH^0(I_{j}) \leq
10^n\omega_n\, \cH^0(I_{j-1})\,\tau^{-k}.
\end{equation}

\medskip

In any case the procedure ends at $j=p$ with a covering of $\cT_{\sfrac{\tau^pr_0}{5}}(S_A)$ which is made of
balls $\{B_{\tau^pr_0}(x_i)\}_{i\in I_p}$ such that $x_i\in S_A$ and
\begin{align}\label{e:cardinality finale}
\cH^0(I_p) &\leq
5^n\tau^{-nq}r_0^{-n}\big(\diam(\Omega)+1\big)^n\big(20^n\,\tau^{-n}\big)^{M}\,
\big(10^n\omega_n\,\tau^{-k}\big)^{p-q-M}\notag\\
& \leq C\, \tau^{-kp}(20^{n} \omega_n)^p\leq C\,\tau^{-p\left(k+\frac{\kappa_0}{2}\right)}
\end{align}
with $C=C(\kappa_0, \delta, r_0, n,\Omega)>0$ and where we used \eqref{e:tau} in the last inequality.
Estimate \eqref{e:stima basica} follows at once
\begin{equation*}
\vert \cT_r(S_A)\vert \leq \cH^0(I_p)\,
|B_{\tau^pr_0}| \stackrel{\eqref{e:cardinality finale}}{\leq}
C\,r^{n-k-\frac{\kappa_0}{2}},
\end{equation*}
for some $C=C(\kappa_0,\delta, r_0, n,\Omega)>0$.

\medskip

\noindent\emph{Step 4: proof of \eqref{e:countable}.}
Let $j_x$ be the smallest index such that \eqref{e:boni in du cess} holds for every $j \geq j_x$, and for every $i\in \N$ set
$$
A_i:=\{x\in\mathcal{S}^0_{r_0, \delta} \;:\; j_x = i\}.
$$
We will prove that $A_i$ is discrete, and hence $\mathcal{S}^0_{r_0, \delta}$ is at most countable. Fix $x\in A_i$. By the choice 
of the parameters applying the Structural Hypothesis (i) it follows that ${\rm{d}}_0(x,4r_0\tau^j)\leq\gamma_0$ for every $j\geq i$. Since $x\in\mathcal{S}^0_{r_0, \delta}$, we can apply Lemma~\ref{l:cattura} and infer that the points $y\in B_{r_0\tau^j}(x)$ satisfying ${\rm{d}}_0(y,4r_0\tau^j)\leq\gamma_0$ are contained in $B_{r_0\tau^{j+1}}(x)$.
Therefore $A_i \cap B_{r_0\tau^j}(x) \subset B_{r_0\tau^{j+1}}(x)$ for every $j \geq i$, which implies that $A_i$ is discrete.
\end{pf}

\subsection{Almgren's stratification principle}\label{s:bianco}
In this section we make the connection to the approach to Almgren's stratification 
principle by White \cite{Wh}. 
Indeed under very natural assumptions the results by White for the time independent 
case follow from ours.

White's stratification criterion in its simplest formulation is based on:
\begin{itemize}
\item[(a$'$)] an upper semi-continuous function $f:\Omega\to[0,\infty)$
defined on a bounded open set $\Omega\subset\R^n$;
\item[(b$'$)] for every $x \in \Omega$ a compact class of conical functions $\cG(x)$
according to the following definition.
\end{itemize}

\begin{definition}\label{d:conical-1}
(1) An upper semi-continuous map $g:\Rn\to[0,\infty)$ is \emph{conical}
if $g(z)=g(0)$ implies that $g$ is positively $0$-homogeneous
with respect to $z$, i.e.,
\[
g(z+\lambda x)=g(z+x)\; \text{ for all $x\in\Rn$ and $\lambda>0$}.
\]
(2) A class $\mathscr{G}$ of conical functions is \emph{compact} if for
all sequences $(g_i)_{i\in\N}\subseteq\mathscr{G}$ there exist
a subsequence $(g_{i_j})_{j\in\N}$
and an element $g\in\mathscr{G}$ such that
\[
\limsup_{j\to \infty} g_{i_j}(y_{i_j})\leq g(y)
\quad \forall \;  y\in \R^n,\;(y_i)_{i\in \N} \subset\R^n \;
\text{ with }\;y_i\to y.
\]
\end{definition}
In particular a conical function is $0$-homogeneous
with respect to $0$.

The stratification theorem by White is then based on the following two
structural hypotheses:
\begin{itemize}
\item[(i$'$)] $g(0)=f(x)$ for all $g\in\mathscr{G}(x)$;
\item[(ii$'$)] for all $r_i\downarrow 0$ there exist a subsequence $r_{i_j}\downarrow 0$
and $g\in\mathscr{G}(x)$ such that
\[
\limsup_{j\to + \infty}f(x+r_{i_j}y_j)\leq g(y)\quad \text{for all $y,y_j\in B_1$ with $y_j\to y$}.
\]
\end{itemize}

By the upper semi-continuity of any conical function $g$, the closed set
\[
S_g:=\{z\in\Rn:\,g(z)=g(0)\}
\]
is in fact the set of the maximum points of $g$.
$S_g$ is called the \emph{spine} of $g$.
Moreover $S_g$ is  the largest vector space that leaves $g$ invariant, i.e.,
\begin{equation}\label{e:sg}
S_g=\{z\in\Rn:\,g(y)=g(z+y)\,\text{ for all }y\in\Rn\}
\end{equation}
(cp.~\cite[Theorem 3.1]{Wh}).
We set $d(x):=\sup\{\dim S_g:g\in\mathscr{G}(x)\}$, and
\[
\Sigma_\ell:=\{x\in \Omega:\,f(x)>0,\,d(x)\leq \ell\}.
\]

The stratification criterion in \cite[Theorem 3.2]{Wh} is
the following.

\begin{thm}[White]\label{t:bianco}
Under the Structural Hypotheses \textup{(i$'$)}, \textup{(ii$'$)},
\begin{gather}
\Sigma_0 \;\text{ is countable};\label{e:gac ipotetico}\\
\dim_{\cH}(\Sigma_\ell) \leq \ell \quad\forall\;\ell\in\{1,\ldots, n\},
\end{gather}
where $\dim_{\cH}$ denotes the Hausdorff dimension.
\end{thm}

The reader who is interested in the application of this criterion to the model cases of area minimizing currents 
and harmonic maps is referred to \cite{Wh}.

\subsubsection{Relation between the structural hypotheses}\label{pa.min bianco}
Theorem~\ref{t:bianco} can be recovered from our Theorem~\ref{t:nostra stratiFICA}
if we assume the following relations between the Structural Hypotheses (i), (ii)
in \S~\ref{sh} and (i$'$), (ii$'$) in \S~\ref{s:bianco}:
\begin{itemize}
\item[(1)] $f = \Theta_0$;
\item[(2)] for every $x \in \Omega$, if
\[
\lim_j  \d_k(x,r_j) =0 \quad\text{for some $(r_j)_{j \in \N} \subset (0, \dist(x,\de\Omega))$,}
\]
then $x \notin \Sigma_{k-1}$.
\end{itemize}
Note that (1) and (2) are always satisfied in the relevant examples
considered in the literature.

To prove that the conclusions of Theorem~\ref{t:bianco}
are implied by Theorem~\ref{t:nostra stratiFICA} it is enough to show that
\begin{equation}\label{e:inclusione2}
\Sigma_\ell \subset \bigcup_{r_0>0} \mathcal{S}^\ell_{r_0}.
\end{equation}
This means that for every $r_0>0$ and
for every $x \in \Sigma_\ell \cap \Omega^{r_0}$ there exists $\delta >0$ such that
\begin{equation}\label{e:inclusione parafrasi}
\d_{\ell+1}(x, r) \geq \delta \quad \forall \; 0<r\leq r_0.
\end{equation}
Assume by contradiction that \eqref{e:inclusione parafrasi} does not hold,
we find $r_0$ and $x$ as above such that
for a sequence $r_j \in (0,r_0]$ we have $\d_{\ell+1}(x, r_j) \downarrow 0$.
Then by \S~\ref{pa.min bianco} (2) $x$ cannot belong to $\Sigma_\ell$.

%
%
\section{Preliminary results on Dir-minimizing Q-valued functions}\label{s:Q}
We follow \cite{DSmemo} for the notation and the terminology, which we
briefly recall in the following subsections.

The space of $Q$-points of $\R^m$ is the subspace of positive atomic measures in $\R^m$ with mass $Q$, i.e.
\[
\Iq(\R^m) := \left\{ \sum_{i=1}^Q \a{p_i} \; : \; p_i \in \R^m \right\}
\]
where $\a{p_i}$ denotes the Dirac delta at $p_i$.
$\Iq$ is endowed with the complete metric $\cG$ given by: for every $T = \sum_i \a{p_i}$ and $S = \sum_i \a{p'_i} \in \Iq(\R^m)$
\begin{gather*}
\cG(T, S) := \min_{\sigma \in \Pe_Q} \left(\sum_{i=1}^Q |p_i - p'_{\sigma(i)}|^2 \right)^{\sfrac12}
\end{gather*}
where $\Pe_Q$ is the symmetric group of $Q$ elements.

A $Q$-valued function is a measurable map $u:\Omega \to \Iq(\R^m)$ from a bounded open set $\Omega\subset \R^n$ (with smooth boundary $\de \Omega$ for simplicity).
It is always possible to find measurable functions $u_i :\Omega \to \R^m$
for $i\in\{1,\ldots, Q\}$ such that $u(x) = \sum_i \a{u_i(x)}$ for a.e.~$x \in \Omega$.
Note that the $u_i$'s are not uniquely determined: nevertheless, we often
use the notation $u = \sum_i \a{u_i}$ meaning an admissible choice
of the functions $u_i$'s has been fixed. We set
\[
|u|(x) := \cG(u(x), Q\a{0}) = \left(\sum_i |u_i(x)|^2\right)^{\sfrac12}.
\]
The definition of the Sobolev space $W^{1,2}(\Omega, \Iq)$ is given in \cite[Definition~0.5]{DSmemo} and leads to the notion of approximate differential $Du = \sum_i \a{Du_i}$ (cp.~\cite[Definitions~1.9 \& 2.6]{DSmemo}.
We set
\begin{gather*}
|Du|(x) := \left(\sum_i |Du_i(x)|^2\right)^{\sfrac12}
\end{gather*}
and say that a function $u \in W^{1,2}(\Omega,\Iq(\R^m))$ is \textit{Dir-minimizing} if
\[
\int_{\Omega} |Du|^2 \leq \int_{\Omega} |Dv|^2 \quad \forall\; v \in W^{1,2}(\Omega),\;v\vert_{\de\Omega} = u\vert_{\de\Omega}
\]
where the last equality is meant in the sense of traces (cp.~\cite[Definition~0.7]{DSmemo}).
By \cite[Theorem~0.9]{DSmemo} Dir-minimizing $Q$-valued functions are locally
H\"older continuous with exponent $\beta=\beta(n,Q)>0$.

In what follows we shall always assume that $u \in W^{1,2}(\Omega,\Iq(\R^m))$ is a nontrivial Dir-minimizing function, 
i.e.~$u\not\equiv Q\a{0}$, with
\begin{equation}\label{e:media 0}
\etaa \circ u := \frac{1}{Q}\sum_{i=1}^Q u_i \equiv 0.
\end{equation}
As explained in \cite[Lemma~3.23]{DSmemo} the mean value condition in \eqref{e:media 0} does not introduce any 
substantial restriction on the space of Dir-minimizing functions.
Moreover, in this case $\Delta_Q$ reduces to the set $\{x\in\Omega:\,u(x)=Q\a{0}\}$.
Note that, if $u\not\equiv Q\a{0}$, then $\Delta_Q \subset \Sing_u$ by \cite[Theorem~0.11]{DSmemo}.

\subsection{Frequency function}

We start by introducing the following quantities: for every $x\in \Omega$ and $s>0$ such that $B_s(x) \subset \Omega$ we set
\begin{gather*}
D_u(x,s) := \int_{B_s(x)}|Du|^2\\
H_u(x,s) := \int_{\de B_s(x)}|u|^2\\
I_u(x,s):=\frac{s\,D_u(x,s)}{H_u(x,s)}.
\end{gather*}
$I_u$ is called the \textit{frequency function} of $u$.
Since $u$ is Dir-minimizing and nontrivial, it holds that $H_u(x,s)>0$ for every $s \in (0,\dist(x,\de \Omega))$ (cp.~\cite[Remark~3.14]{DSmemo}), from which $I_u$ is well-defined.

We recall that the functions $s\mapsto D_u(x,s)$, $s\mapsto H_u(x,s)$, and $s\mapsto I_u(x,s)$ are absolutely continuous on $(0,\dist(x,\de \Omega))$.
Similarly for fixed $s \in (0,\dist(x,\de \Omega))$ one can prove the continuity of $x \mapsto D_u(x,s)$,  $x \mapsto H_u(x,s)$ and $x \mapsto I_u(x,s)$ for $x \in \{ y : \dist(y,\de\Omega) >s\}$. 
The former follows by the absolute continuity of Lebesgue integral; while for the remaining two it suffices the following estimate:
\begin{align}\label{e:conto H}
\left\vert\sqrt{H_u(x,s)} - \sqrt{H_u(y,s)} \right\vert& \leq
\left( \int_{\de B_s(y)} ||u|(z) - |u|(z+x-y)|^2\,dz
\right)^{\frac{1}{2}}\notag\\
&\leq |x-y|\left(\int_{\de B_s(y)} \int_0^1|\nabla |u| (z+t\,(x-y))|^2
\,dt\,dz\right)^{\frac{1}{2}}\notag\\
&\leq |x-y|
\left(\int_{B_{s+|x-y|}(y)\setminus B_{s-|x-y|}(y)}|Du|^2\right)^{\frac12}
\end{align}
where we use the fact that $|u| \in W^{1,2}(\Omega)$ with $|\nabla |u|| \leq |Du|$ (cp.~\cite[Definition~0.5]{DSmemo}).

The following monotonicity formula discovered by Almgren in \cite{Alm} is the main estimate about Dir-minimizing functions (cp.~\cite[Theorem~3.15 \& (3.48)]{DSmemo}): for all $0 \leq r_1 \leq r_2 < \dist(x, \de \Omega)$ it holds
\begin{multline}\label{e:Almgren formulozzo}
I_u(x,r_2) - I_u(x,r_1)\\
= \int_{r_1}^{r_2} \frac{t}{H_u(t)}
\Big(\int_{\de B_t(x)}|\de_\nu u|^2 \,\int_{\de B_t(x)}|u|^2
- \Big( \int_{\de B_t(x)}\langle \de_\nu u, u\rangle\Big)^2 \Big)\,dt.
\end{multline}

We finally recall that from \cite[Corollary~3.18]{DSmemo} we also deduce that
\begin{equation}\label{e.asintotica}
H_u(z,r) = O(r^{n+2\,I_u(z,0^+)-1})
\end{equation}
where $I_u(z,0^+) = \lim_{r\downarrow 0} I_u(z,r)$.

\subsection{Compactness}\label{pa.cpt}
From \cite[Proposition~2.11 \& Theorem~3.20]{DSmemo}, if $(u_j)_{j \in \N}$ is a sequence of Dir-minimizinig functions in $\Omega$ such that
\[
\sup_j\|u_j\|_{L^2(\Omega)} + \sup_j\int_\Omega|Du_j|^2 < + \infty,
\]
then there exists $u \in W^{1,2}(\Omega, \Iq)$ such that $u$ is Dir-minimizing, and
up to passing to a subsequence (not relabeled in the sequel) $\cG(u_j,u)\to0$ in $L^2(\Omega)$, and for every $\Omega'\subset\subset \Omega$
\[
\|\cG(u_j, u)\|_{L^\infty(\Omega')} \to 0
\quad\text{and}\quad
\int_{\Omega'}|Du_j|^2 \to \int_{\Omega'}|Du|^2.
\]
In particular this implies that $(|Du_j|^2)_{j\in \N}$ are equi-integrable in $\Omega'$,
and
\begin{equation}\label{e.continuity frequency}
\lim_{j\to +\infty}I_{u_j}(x,s) = I_u(x,s) \quad \forall\;x\in \Omega, \; \forall\;0<2 s< \dist(x,\de \Omega).
\end{equation}

\subsection{Homogeneous $Q$-valued functions}\label{s.bu}
We discuss next some properties of the class of homogeneous $Q$-valued functions: 
$w \in W^{1,2}_\loc(\R^n, \Iq(\R^m))$ satisfying
\begin{itemizeb}
\item[(1)] $w$ is locally Dir-minimizing with $\etaa \circ w \equiv 0$,
\item[(2)] $w$ is $\alpha$-homogeneous, in the sense that
\[
w(x) = |x|^\alpha \, w\left(\frac{x}{|x|}\right) \quad
\forall \, x \in \R^n\setminus\{0\},
\]
for some $\alpha \in (0, \Lambda_0]$, where $\Lambda_0$ is a constant to be specified later.
\end{itemizeb}
We denote this class by $\mathcal{H}_{\Lambda_0}$.
Note that $I_w(x,0^+) = 0$ if $w(x) \neq Q\a{0}$.
The following lemma is an elementary consequence of the definitions.

\begin{lem}\label{l.tangenti}
Let $w \in \cH_{\Lambda_0}$.
Then $I_w(\cdot,0^+)$ is conical in the sense of Definition~\ref{d:conical-1} (1).
\end{lem}

\begin{pf}
Firstly $I_w(\cdot,0^+)$ is upper semi-continuous. Indeed since $w$ is Dir-minimizing, we can use \eqref{e:Almgren formulozzo} and deduce that $I_w(x,0^+)=\inf_{s>0} I_w(x,s)$, i.e.~$I_w(\cdot,0^+)$ is the infimum of continuous (by \eqref{e:conto H}) functions $x\mapsto I_w(x,s)$ and hence upper semi-continuous.

We need only to show that $I_w(\cdot,0^+)$ is $0$-homogeneous at every point $z$ such that $I_w(z,0^+) = I_w(0,0^+)$. We can assume without loss of generality that $w$ is nontrivial, i.e.~$w \not\equiv Q\a{0}$.
We start noticing that if $I_w(z,0^+) = I_w(0,0^+)$ then
\[
I_w(z,0^+) = I_w(0,0^+) = I_w(0, 1) >0
\]
where in the last equality we used the homogeneity of $w$.
Therefore in particular $w(z) = Q\a{0}$. Next we show that $I_w(z,r) = I_w(0,0^+)$ for all $r>0$.
By a simple estimate we get
\begin{align}\label{e.dall'alto}
I_w(z,r) &= \frac{r\,D_w(z,r)}{H_w(z,r)}
\leq I_w(0,r+|z|)\;\frac{H_w(0,r+|z|)}{H_w(0,r)}
\;\frac{H_w(0,r)}{H_w(z,r)}.
\end{align}
Since $w$ is homogeneous with respect to the origin and the frequency of $w$ at $0$ is exactly $\alpha$ (cp.~\cite[Corollary~3.16]{DSmemo}), we have also
\begin{gather*}
H_w(0,r) = H_w(0,1)\,r^{n+2\,\alpha-1}\\
D_w(0,r) = D_w(0,1)\,r^{n+2\,\alpha-2}.
\end{gather*}
In particular
\begin{gather*}
I_w(0,r+|z|) = \alpha = I_w(0,0^+) = I_w(z,0^+)\\
\frac{H_w(0,r+|z|)}{H_w(0,r)} \to 1 \quad\text{as} \quad r\uparrow+\infty.
\end{gather*}
For what concerns the third factor in \eqref{e.dall'alto}
\begin{align}
\frac{H_w(0,r)}{H_w(z,r)}
=  1+ \frac{H_w(0,r)-H_w(z,r)}{H_w(z,r)}
\end{align}
and from \eqref{e.asintotica} and \eqref{e:conto H} we infer that
\begin{align}\label{e:H centri diversi}
|H_w(0,r)-H_w&(z,r)| = \big(\sqrt{H_w(0,r)} + \sqrt{H_w(z,r)}\big)
|\sqrt{H_w(0,r)} - \sqrt{H_w(z,r)}|\notag\\
&\leq C\, r^{\frac{n+2\,I_u(0,0^+)-1}{2}}\,
|z|
\big(D_w(0, {r+|z|}) - D_w(0, {r-|z|})\big)^{\frac12}\notag\\
&\leq C\,|z|\, r^{\frac{n+2\,I_u(0,0^+)-1}{2}}\,
\big((r+|z|)^{n+2\,\alpha -2} - (r-|z|)^{n+2\,\alpha -2}\big)^{\frac12}\notag\\
&\leq C\,|z|^{\frac32}\,r^{n+2\,\alpha - 2}.
\end{align}
This in turn implies
\[
\frac{H_w(0,r)}{H_w(z,r)} \to 1 \quad \text{as }\;r\uparrow+\infty
\]
and from \eqref{e.dall'alto}
\begin{equation*}
\lim_{r \to + \infty} I_w(z, r) \leq \lim_{r\downarrow 0^+} I_w(z, r),
\end{equation*}
i.e.~by Almgren's monotonicity estimate \eqref{e:Almgren formulozzo} we infer that $I_w(z, r) = I_w(z,0^+)$ for all $r>0$. 
As a consequence (cp.~\cite[Corollary~3.16]{DSmemo}) $w$ is $\alpha$-homogeneous at $z$ which straightforwardly implies that $I_w(\cdot,0^+)$ is $0$-homogeneous at $z$.
\end{pf}

We can then define the spine of a homogeneous $Q$-valued function $w \in \cH_{\Lambda_0}$:
\[
S_w := \{ x \in \R^n \;:\; I_w(x,0^+) = I_w(0,0^+) \}.
\]
By the proof of Lemma~\ref{l.tangenti} it follows that $w$ is $\alpha$-homogeneous at every point $x \in S_w$. Similarly it is simple to verify that $S_w$ is the largest vector space which leaves $w$ invariant, as well as $I_w(\cdot,0^+)$:
\begin{equation}\label{e.invariante}
S_w = \big\{ z \in \R^n \; :\; w(y) = w(z + y) \quad \forall\; y \in \R^n \big\}.
\end{equation}
Indeed it is enough to prove that every $z \in S_w$ leaves $w$ invariant (the other inclusion is obvious). To show this, note that by the $\alpha$-homogeneity of $w$ at $z$ and $0$ it follows that for every $y \in \R^n$
\begin{align*}
w(y) & = w\left(z + y - z\right) = 2^\alpha \, w\left(z + \frac{y - z}{2}\right) = 2^\alpha \, w\left(\frac{y + z}{2}\right)\\
& = w\left(z + y\right).
\end{align*}

We denote by $\cC_k$ for $k \in\{ 0, \ldots, n\}$ the set of $k$-invariant homogeneous $Q$-functions
\begin{equation}\label{e:CkQ}
\cC_k := \{w \in \cH_{\Lambda_0}\;:\;\dim (S_w) \geq k \}.
\end{equation}
Note that $\cC_n = \cC_{n-1} = \{Q\a{0}\}$, i.e.~these sets are singleton consisting of the constant function $w \equiv Q\a{0}$.
For $\cC_n$ this is follows straightforwardly from the definition and \eqref{e.invariante}. While for $\cC_{n-1}$ one can argue via the cylindrical blowup in \cite[Lemma~3.24]{DSmemo}. Assume without loss of generality that
\[
w \in \cC_{n-1},  \quad w  \not\equiv Q\a{0} \quad \text{and}\quad S_w = \R^{n-1} \times \{0\}.
\]
Then by the invariance of $w$ along $S_w$ it follows that $w$ is a function of one variable.
By \cite[Lemma~3.24]{DSmemo} it follows that $\tilde w: \R \to \Iq(\R^m)$ is locally Dir-minimizing and 
\[
\tilde w\not\equiv Q\a{0}, \quad \etaa\circ \tilde w \equiv 0.
\]
This is clearly a contradiction because the only Dir-minimizing function of one variable are non-intersecting linear functions (cp.~\cite[3.6.2]{DSmemo}).

Finally, a simple consequence of \eqref{e.invariante} is that $\{ w \vert_{B_1} : w \in \cC_k\}$ is a closed subset of $L^2(B_1, \Iq(\R^m))$.

\begin{lem}\label{l.chiusura hom}
Let $(w_j)_{j \in \N} \subset \cC_k$ and $w \in W^{1,2}_\loc(\Omega,\Iq(\R^m))$ be such that $w_j \to w$ in $L^2_\loc(\R^n, \Iq(\R^m))$. Then $w \in \cC_k$.
\end{lem}

\begin{pf}
Let $\alpha_j$ be the homogeneity exponent of $w_j$.
Since for Dir-minimizing $\alpha$-homogeneous $Q$-valued functions $w$ it holds that $D_w(1) = \alpha\, H_w(1)$, we deduce from $\alpha_j \leq \Lambda_0$ and $w_j \to w$ that the functions $w_j$ have equi-bounded energies in any compact set of $\R^n$. Therefore by the compactness in \S~\ref{pa.cpt} it follows that $w_j \to w$ locally uniformly and $w \in \cH_{\Lambda_0}$.

For every $j \in \N$ let now $V_j$ be a $k$-dimensional linear subspace of $\R^n$ contained in $S_{w_j}$. By the compactness of the Grassmannian $\textup{Gr}(k,n)$, we can assume that up to passing to a subsequence (not relabeled) $V_j$ converges to a $k$-dimensional subspace $V$. Using the uniform convergence of $w_j$ to $w$ we then conclude that for every $z \in V$ and $y \in \R^n$
\[
w(z+y) = \lim_{j}w_j(z_j + y) = \lim_j w_j(y) = w(y)
\]
where $z_j \in V_j$ is any sequence such that $z_j \to z$.
This shows that $V \subset S_w$, thus implying that $\dim(S_w) \geq k$.
\end{pf}

\subsection{Blowups}
Let $u$ be a Dir-minimizing $Q$-valued function, $\etaa \circ u \equiv 0$ and
$u \not\equiv Q\a{0}$. Fix any $r_0>0$.
For every $y\in \Delta_Q\cap \Omega^{r_0}$, i.e.~for every $y$ such that $u(y) = Q\a{0}$
and $\dist(y, \de \Omega) \geq 2 r_0$, we define the rescaled functions of $u$ at $y$ as
\[
u_{y,s}(x) := \frac{s^{\frac{m-2}{2}}u(y+sx)}{D^{\sfrac12}_u(y,s)} \qquad \forall\; 0 < s < r_0, \; \forall \; x \in B_{\frac{r_0}{s}}(0). 
\]
From \cite[Theorem~3.20]{DSmemo} we deduce that for every $s_k \downarrow 0$ there exists a subsequence 
$s_k' \downarrow 0$ such that $u_{y,s'_k}$ converges locally uniformly in $\R^n$ to a function 
$w:\R^n \to \Iq(\R^m)$ such that $w \in  \cH_{\Lambda_0}$ with 
\begin{equation}\label{e:Lambda_0}
\Lambda_0=\Lambda_0(r_0):= \frac{r_0 \int_\Omega |Du|^2}{\min_{x \in \Omega^{r_0}} H_u(x,r_0)}.
\end{equation}
Note that $\min_{x \in \Omega^{r_0}} H_u(x,r_0) >0$. Indeed, by the continuity of $x \mapsto H_u(x,r_0)$ and the closure 
of $\Omega^{r_0}$, the minimum is achieved and cannot be $0$ because of the condition $u \not\equiv 0$. 
In particular, $\Lambda_0 \in \R$.

%
%

\section{Stratification for Dir-minimizing $Q$-valued functions}\label{s:dirminstrat}
In this section we apply Theorems~\ref{t:tub-neig}, \ref{t:nostra stratiFICA} and \ref{t:minchioschi} to the case of Almgren's Dir-minimizing $Q$-valued functions.
Keeping the notation $\Omega^s$ and $U$ as in \S~\ref{pa:1.1}, we set
\begin{itemizeb}
\item[(1)] $\Theta_s: \Omega^s \to [0,+\infty)$ given by
\[
\Theta_0(x):= \lim_{r\downarrow 0^+} I_u(x,r) \quad \text{and}\quad
\Theta_s(x):= I_u(x,s) \quad \text{for }\; s>0, x\in\Omega^s\;,
\]
\item[(2)] for every $k\in\{0, \ldots, n\}$, $\d_k : U^{} \to [0,+\infty)$ is given by
\[
\d_k(x,s) := \min \Big\{ \|\cG(u_{x,s}, w)\|_{L^2(\de B_1)} : w \in \cC_k \Big\}.
\]
Note that since $\{ w \vert_{B_1} : w \in \cC_k\}$ is a closed subset of $L^2(B_1)$ the minimum in the definition of $d_k$ is achieved.
\end{itemizeb}

It follows from Almgren's monotonicity formula \eqref{e:Almgren formulozzo} that conditions (a) and (b) of \S~\ref{pa:1.1} are satisfied.

\medskip

We verify next that the Structural Hypotheses in \S~\ref{sh} are fulfilled. 
For simplicity we write the corresponding statements for fixed $r_0$. The corresponding $\Lambda_0>0$ 
is defined as in \eqref{e:Lambda_0} above. Therefore, the sets $\mathcal{H}_{\Lambda_0}$ and $\cC_k$, 
introduced respectively in \S~\ref{s.bu} and \eqref{e:CkQ}, are defined in terms of $\Lambda_0=\Lambda_0(r_0)$.

\begin{lem}\label{l:Dir:(i)}
For every $\eps_1>0$ there exist $0<\lambda_1(\eps_1),\,\eta_1(\eps_1) <\sfrac{1}{4}$ such that, for all $(x,s)\in U$ with $x\in\Omega^{r_0}$ and $s<r_0$, it holds
\[
I_u(x, s) - I_u(x, \lambda_1 s)
\leq\eta_1 \quad\Longrightarrow\quad
\exists \; w \in \cC_0 \;:\; \|\cG(u_{x,s},w)\|_{L^2(\de B_1)} \leq \eps_1.
\]
\end{lem}

\begin{pf}
We argue by contradiction and assume there exist points $(x_j, s_j)$ with $x_j\in\Omega^{r_0}$ and $s_j<r_0$ such that
\begin{equation*}
I_u(x_j, s_j) - I_u(x_j, \tfrac{s_j}{2^j}) \leq 2^{-j} \quad\text{and}\quad
\|\cG(u_{x_j,s_j},w)\|_{L^2(\de B_1)} \geq \eps_1 \quad \forall\; w\in \cC_0
\end{equation*}
or equivalently, setting $u_j:= u_{x_j,s_j}$, 
\begin{equation}\label{e:contra2}
I_{u_j}(0,1)-I_{u_j}(0,2^{-j})\leq 2^{-j}
\quad\text{and}\quad
\|\cG(u_{j},w)\|_{L^2(\de B_1)} \geq \eps_1 \quad \forall\; w\in \cC_0.
\end{equation}
From \cite[Corollary~3.18]{DSmemo} it follows that
\begin{equation}\label{e:equi bounded}
\sup_j D_{u_j}(0,2) \leq 2^{n-2+2\, I_{u_j}(0,2)}\,\frac{I_{u_j}(0,2)}{I_{u_j}(0,1)} \leq C
\end{equation}
where $C = C(\Lambda_0)$ because $I_{u_j}(0,2) \leq \Lambda_0$ by definition of $\Lambda_0$.
We can then use the compactness for Dir-minimizing functions in \S~\ref{s.bu} to infer the existence of a Dir-minimizing $w$ such that (up to subsequences) $u_j \to w$ locally strongly in $W^{1,2}(B_2)$ and uniformly.
We then can pass into the limit in \eqref{e:Almgren formulozzo} and using \eqref{e:contra2} we obtain
\[
\int_{0}^{1} \frac{t}{H_{w}(t)}
\left(\int_{\de B_t}|\de_\nu w|^2 \,\int_{\de B_t}|w|^2
- \left( \int_{\de B_t}\langle \de_\nu w, w\rangle\right)^2 \right)\,dt=0.
\]
This implies that $w$ is $\alpha$-homogeneous (cp.~\cite[Corollary~3.16]{DSmemo}) with
\[
\alpha = \lim_{j} I_{u_j}(0,1) \leq \Lambda_0
\]
because of \S~\ref{pa.cpt}.
This contradicts $\|\cG(u_{j},w)\|_{L^2(\de B_1)} \geq \eps_1$ for all $w\in \cC_0$ in \eqref{e:contra2} and proves the lemma.
\end{pf}

\begin{rmk}
Using the regularity theory of Dir-minimizing functions proven in \cite{DSmemo} it is in fact possible to prove a stronger claim then Lemma~\ref{l:Dir:(i)}, namely that for every $\eps_1>0$ there exists $0<\eta_1(\eps_1) <\sfrac{1}{4}$ such that for all $(x,s)\in U$ with $x\in\Omega^{r_0}$ and $s<r_0$
\begin{equation}\label{e.cosa inutile}
I_u(x, s) - I_u(x, \sfrac{s}{2})
\leq\eta_1 \quad\Longrightarrow\quad
\exists \; w \in \cC_0 \;:\; \|\cG(u_{x,s},w)\|_{L^2(\de B_1)} \leq \eps_1.
\end{equation}
Since \eqref{e.cosa inutile} is not needed in the sequel, we leave the details of the proof to the reader.
\end{rmk}

For what concerns (ii) we argue similarly using a rigidity property of homogeneous Dir-minimizing functions.

\begin{lem}\label{l:Dir:(ii)}
For every $0<\eps_2, \tau <1$ there exists $0<\eta_2(\eps_2,\tau)\leq \eps_2$
such that if $(x,5s)\in U$, with $x\in\Omega^{r_0}$ and $5s<r_0$, $\d_k(x,4s) \leq \eta_2$ 
and $\d_{k+1}(x,4s) \geq \eps_2$ for some $k \in \{0, \ldots, n-1\}$ then there exists a $k$-dimensional affine space $V$ such that
\[
\d_0(y,4s)> \eta_2 \quad \forall\; y \in B_{s}(x)\setminus
\cT_{\tau s}(V).
\]
\end{lem}

\begin{pf}
We prove the statement for $V= S_w$ with $w\in \cC_k$ such that $\d_k(x,4s)=\|\cG(u, w)\|_{L^2(\de B_{4s}(x))}$.
We argue by contradiction. Reasoning as above with the rescalings of $u$ (eventually composing with a rotation of the domain to achieve (4) below for a fixed space $V$), we find a sequence of functions $u_j \in W^{1,2}(B_5, \Iq(\R^k)$ such that
\begin{itemizeb}
\item[(1)] $\sup_j D_{u_j}(0,5) < + \infty$;
\item[(2)] there exists $w_j \in \cC_k$ such that $\|\cG(u_j, w_j)\|_{L^\infty(B_4)} \downarrow 0$;
\item[(3)] $\|\cG(u_j, w)\|_{L^2(B_4)} \geq \eps_2$ for every $w \in \cC_{k+1}$;
\item[(4)] there exists $y_j \in B_1\setminus \cT_\tau(V)$ such that $\d_0(y_j,4) \downarrow 0$ and $V=S_{w_j}$ is the $k$-dimensional spine of $w_j$ (note that by (2) \& (3) the dimension of the spine of $w_j$ cannot be higher than $k$). 
\end{itemizeb}
Possibly passing to subsequences (as usual not relabeled), we can assume that $u_j \to w$, $w_j \to w$ locally in $L^{2}(\R^n,\Iq(\R^m))$ and $y_j \to y$ for some $w \in W^{1,2}_{\textup{loc}}(\R^n, \Iq(\R^m))$ and $y \in \bar B_1\setminus \cT_\tau(V)$.
By Lemma~\ref{l.chiusura hom} we deduce that $w \in \cC_k$ with $S_w \supset V$; since by (3) $w\not\in \cC_{k+1}$, we conclude $S_w = V$.

It follows from (4) that $w_{y,s} = w_{y,1}$ for every $s\in (0,1]$.
Indeed there exist $z_j \in \cC_0$ such that $\|\cG((u_j)_{y_j,1}, z_j)\|_{L^2(\de B_4)} \downarrow 0$ and by continuity $(u_{j})_{y_j, 1} \to w_{y,1} \in \cC_0$.
In particular $w(y)=0$ and by the upper semi-continuity of $x\mapsto I_w(x,0^+)$ we deduce also that $I_w(y,0^+) = I_w(0,0^+)$, i.e.~$y \in S_w$ which is the desired contradiction.
\end{pf}

We can then infer that Theorem~\ref{t:tub-neig} holds for $Q$-valued functions. 

\begin{thm}\label{t:Q neig}
Let $u :\Omega \to \Iq(\R^m)$ be a nontrivial Dir-minimizing function with average $\etaa\circ u \equiv 0$.

For every $0<\kappa_0, \delta <1$  and $r_0 >0$, there exists $C=C(\kappa_0, \delta, r_0, n)>0$ such that
\begin{gather*}
|\cT_r(\Delta_Q\cap \mathcal{S}^{k}_{r, r_0, \delta})|\leq C\, r^{n-k-\kappa_0} \quad \forall \; 
k \in \{1, \ldots, n-1\},\\
\text{and}\quad \mathcal{S}^0_{r_0, \delta} \text{ is countable}.
\end{gather*}
\end{thm}

In particular, Theorem~\ref{t:nostra stratiFICA} applies and we conclude that $\dim_\cH(\mathcal{S}^k_{r_0}) \leq k$ 
and that $\mathcal{S}^0_{r_0}$ is at most countable.
We shall improve upon the latter estimate on the stratum $\mathcal{S}^{n-1}_{r_0}$ in the next section.

\subsection{Minkowski dimension}
We can actually give an estimate on the Minkowski dimension of the set of maximal multiplicity points $\Delta_Q$
by means of Theorem~\ref{t:minchioschi}. An $\eps$-regularity result is the key tool to prove this.

\begin{proposition}\label{p:eps reg}
There exists a constant $\delta_0 = \delta_0(r_0)>0$ such that
\begin{equation}\label{e.eps reg minchione}
\mathcal{S}^{n-1}_{r}=\mathcal{S}^{n-2}_{r} = \mathcal{S}^{n-2}_{r, \delta_0} \quad \forall \; r \in (0,r_0).
\end{equation}
\end{proposition}

\begin{pf}
The first equality is an easy consequence of $\cC_n=\cC_{n-1}=\{Q\a{0}\}$ that gives $\d_n\equiv\d_{n-1}$.

Set $\delta_0 := (\Lambda_0 + 1)^{-\sfrac{1}{2}}$, we show that 
$\mathcal{S}^{n-2}_{r,\delta} \subset \mathcal{S}^{n-2}_{r, \delta_0}$ for every $\delta \in (0, \delta_0)$.
Assume by contradiction that there exists $x \in \mathcal{S}^{n-2}_{r,\delta} \setminus \mathcal{S}^{n-2}_{r, \delta_0}$ 
for some $\delta$ as above. From $\cC_{n-1} = \{ Q\a{0}\}$ we deduce the existence of $s \in (0,r)$ such that
\[
0<\delta \leq \| u_{x,s} \|_{L^2(\de B_1)} < \delta_0.
\]
In particular, 
the condition $\int_{B_1}|D u_{x,s}|^2 = 1$ gives
\[
I_{u_{x,s}}(0,1) = \frac{\int_{B_1}|D u_{x,s}|^2}{\int_{\de B_1}|u_{x,s}|^2} \geq \frac{1}{\delta_0^2} > \Lambda_0.
\]
By recalling that $I_u(x,s) = I_{u_{x,s}}(0,1)$, the desired contradiction follows from Almgren's monotonicity formula 
\eqref{e:Almgren formulozzo} and the very definition of $\Lambda_0$ in \eqref{e:Lambda_0}.
\end{pf}

In particular Theorem~\ref{t:A} follows from Theorem~\ref{t:minchioschi}.

\begin{pf}[of Theorem~\ref{t:A}]
It is a direct consequence of Proposition~\ref{p:eps reg} and Theorem~\ref{t:minchioschi}. 
Given $u :\Omega \to \Iq(\R^m)$ a nontrivial Dir-minimizing function (i.e.~$\Delta_Q \neq \Omega$),
we can consider the function 
\[
v(x) := \sum_i \a{u_i(x) - \etaa\circ u(x)}.
\]
Then by \cite[Lemma~3.23]{DSmemo} $v$ is Dir-minimizing with $\etaa\circ v \equiv 0$.
Moreover, the set of $Q$-multiplicity points of $u$ in $\Omega^{r_0}$
corresponds to the set $\mathcal{S}^{n-2}_{r_0}$ for the function $v$ and the conclusion follows
straightforwardly.
\end{pf}

\subsection{Almgren's stratification}
In this section we show that Theorem~\ref{t:bianco} applies in the case of $Q$-valued functions, as well. 
In particular, this implies that the singular strata for Dir-minimizing $Q$-valued functions 
can also be characterized by the spines of the blowup maps, thus leading to the proof of Theorem~~\ref{t:B} 
in the introduction.

By following the notation in \S~\ref{pa.min bianco} (1), we set
\[
f(x) := I_u(x,0^+) \quad \forall\; x \in \Omega. 
\]
For every $x \in \Omega$ such that $f(x) =0$ (or, equivalently, $u(x) \neq Q\a{0}$) we define $\mathscr{G}(x)$ 
to be the singleton made of the constant function $0$, i.e.~$\mathscr{G}(x)= \{Q \a{0}\}$; otherwise
\begin{equation}\label{e:sticazziG}
\mathscr{G}(x) := \big\{I_w(\cdot, 0^+) : w \in W^{1,2}_{\textup{loc}}(\R^n, \Iq(\R^m))\; \text{blowup of } u \text{ at } x \big\}.
\end{equation}
As explained in \S~\ref{s.bu} $\mathscr{G}(x)$ is never empty because there always exist (possibly non-unique) blowup of $u$ at any multiplicity $Q$ point.

Since every blowup of $u$ is a nontrivial homogeneous Dir-minimizing function, it follows from Lemma~\ref{l.tangenti} that every function $g \in \mathscr{G}(x)$ is conical in the sense of Definition~\ref{d:conical-1} (1). We need then to show the following.

\begin{lem}\label{l.cpt}
For every $x \in \Omega$ the class $\mathscr{G}(x)$ is compact in the sense of Definition~\ref{d:conical-1} (2).
\end{lem}

\begin{pf}
If $x$ is not a multiplicity $Q$ point, then there is nothing to prove.
Otherwise consider a sequence of maps $g_j= I_{w_j}(\cdot, 0^+) \in \mathscr{G}(x)$, with $w_j$ blowup of $u$ at $x$.
By \S~\ref{s.bu} $w_j$ is Dir-minimizing $\alpha$-homogeneous with $\alpha = I_u(x,0^+)$ and $D_{w_j}(1) = 1$.
Then by the compactness in \S~\ref{pa.cpt}, there exists $w$ such that $w_j \to w$ locally in $L^{2}$ up to subsequences (not relabeled) with $D_w(1)=1$.
By a simple diagonal argument it follows that $w$ is as well a blowup of $u$ at $x$, i.e.~$g=I_w(\cdot, 0^+) \in \cG(x)$.
For every $y_j \in B_1$ with $y_j \to y \in B_1$ and for every $s>0$, we then deduce
\begin{align*}
\limsup_{j\uparrow+\infty} g_j(y_j) &
\leq \limsup_{j\uparrow+\infty} I_{w_j}(y_j,s)\notag\\
& = \limsup_{j\uparrow+\infty} \left(
\frac{s\,D_{w_j}(y, s)}{H_{w_j}(y,s)}\,\frac{D_{w_j}(y_j,s)}{D_{w_j}(y,s)}\,
\frac{H_{w_j}(y,s)}{H_{w_j}(y_j,s)}\right)\\
&= I_w(y,s)
\end{align*}
where we used
\begin{itemize}
\item[-] the monotonicity of $I_{w_j}(y_j, \cdot)$ in the first line,
\item[-] the continuity of $x \mapsto D_{w_j}(x,s)$ and $x\mapsto H_{w_j}(x,s)$,
\item[-] the convergence of the frequency functions $I_{w_j}(y,s) \to I_w(y,s)$ (cp.~\ref{pa.cpt}).
\end{itemize}
Sending $s$ to $0$ provides the conclusion.
\end{pf}

Finally we prove that the Structural Hypothesis (ii$'$) of White's theorem (cp. \S~\ref{s:bianco}) 
holds as well:
\begin{align*}
\limsup_{j\uparrow+\infty} f(x+r_{i_j} y_j) &=
\limsup_{j\uparrow+\infty} I_u(x+r_{i_j} y_j,0^+)\\
&\leq\limsup_{j\uparrow+\infty} I_u(x+r_{i_j} y_j,r_{i_j}s)\notag\\
&=\limsup_{j\uparrow+\infty}  I_{u_{x,r_{i_j}}}(y_j,s)
= I_{w}(y,s)
\end{align*}
where we used the strong convergence of the frequency of \S~\ref{pa.cpt}.

In particular, Theorem~\ref{t:bianco} holds true, which in turn leads to the proof of Theorem~\ref{t:B}
by a simple induction argument.

\subsection{Stratification Theorem~\ref{t:B}}\label{ss:B}
We define now the singular strata $\Sing^k_u$
for a Dir-minimizing multiple valued function $u:\Omega \to \Iq(\R^m)$. 
Consider any point $x_0 \in \Sing_u$, and let
\[
u(x_0) = \sum_{i=1}^J \kappa_i \, \a{p_i}
\]
with $\kappa_i \in \N\setminus\{0\}$ such that $\sum_{i=1}^J\kappa_i = Q$
and $p_i \neq p_j$ for $i\neq j$.
Then by the uniform continuity of $u$ there exist $r>0$ and
Dir-minimizing multiple valued functions
$u_i:B_r(x_0) \to \I{\kappa_i}(\R^m)$ for $i\in\{1, \ldots, J\}$ such that
\[
u\vert_{B_r(x_0)} = \sum_{i=1}^J \a{u_i},
\]
where by a little abuse of notation the last equality is meant
in the sense $u(x) = \sum_{i}u_i(x)$ as measures.
For every $i\in\{1, \ldots, J\}$ let $v_i:B_r(x_0) \to \I{\kappa_i}(\R^m)$ be
given by 
\[
v_i (x) := \sum_{l=1}^{\kappa_i} \a{(u_i(x))_l -  \etaa\circ u_i(x)}.
\]
Then we say that a point $x_0 \in \Sing_u$
belongs to $\Sing^k_u$, $k\in\{0,\ldots, n\}$,
if the spine of every blowup of $v_i$ at $x_0$, for every $i\in\{1, \ldots, J\}$,
is at most $k$-dimensional.

We can then prove Theorem~\ref{t:B} by a simple induction argument on the number of values $Q$.

\begin{pf}[of Theorem~\ref{t:B}]
Clearly if $Q=1$ there is nothing to prove because every harmonic
function is regular and $\Sing_u = \emptyset$.
Now assume we have proven the theorem for every $Q^* < Q$ and we prove it for $Q$.

We can assume without loss of generality that $\Delta_Q \neq \Omega$.
Then, as noticed, $\Delta_Q=\Sing_u \cap \Delta_Q$ by \cite[Theorem~0.11]{DSmemo}.
Moreover $\Sing^k_u \cap \Delta_Q = \Sigma_k$,
where $\Sigma_k$ is that of Theorem~\ref{t:bianco}.
Indeed $x_0 \in \Sigma_k$ if and only if
the maximal dimension of the spine of any $g\in \mathcal{G}(x_0)$ is at most $k$.
By \eqref{e:sticazziG} $g \in \mathcal{G}(x_0)$ if and only if
$g = I_w(\cdot, 0^+)$ for some blowup $w$ of $u$ at $x_0$.
Hence by \eqref{e.invariante} $x_0 \in \Sigma_k$ if and only if
the dimension of the spines of the blowups of $u$ at $x_0$ is at most $k$.
Note that $\Sing^{n-2}_u \cap \Delta_Q = \Delta_Q$ since $\cC_n = \cC_{n-1} = \{Q\a{0}\}$ 
(we use here the notation in \S~\ref{s.bu}) and $u$ is not trivial.
Therefore we deduce that
\begin{gather*}
\Sing^0_u \cap \Delta_Q \quad \text{is countable}\\
\dim_\cH (\Sing^k_u \cap \Delta_Q) \leq k \quad \forall\; k\in\{1, \ldots, n -2\}.
\end{gather*}

Next we consider the relatively open set $\Omega \setminus \Delta_Q$
(recall that both $\Sing_u $ and $\Delta_Q$ are relatively closed sets).
Thanks to the continuity of $u$ we can find a cover of 
$\Omega \setminus (\Sing_u \cap \Delta_Q)$ made of
countably many open balls $B_i \subset \Omega \setminus (\Sing_u \cap \Delta_Q)$
such that $u\vert_{B_i} = \a{u^1_i} + \a{u^2_i}$ with $u^1_i$ and $u^2_i$
Dir-minimizing multiple valued functions taking strictly less than $Q$ values.
Since $\Sing^k_u \cap B_i = \Sing^k_{u^1_i} \cup \Sing^k_{u^2_i}$ by the very definition,
using the inductive hypotheses for $u^1_i$ and $u^2_i$ we deduce that
\begin{gather*}
\Sing^0_u \cap B_i \quad \text{is countable}\\
\Sing^{n-2}_u \cap B_i = \Sing_u \cap B_i\\
\dim_\cH (\Sing^k_u \cap B_i) \leq k \quad \forall\; k\in\{1, \ldots, n-2\},
\end{gather*}
thus leading to \eqref{e:B1} and \eqref{e:B2}.
\end{pf}

%
%

\section{Applications to generalized submanifolds}\label{s:vrfld}

In the present section we apply the abstract stratification results in \S~\ref{s:abstract}
to integral varifolds with mean curvature in $L^\infty$ and to almost minimizers in
codimension one (both frameworks are not covered by the results in \cite{ChNa2} although 
they can be considered as slight variants of those).
This case is relevant in several variational problems (see the examples in \cite[\S~4]{Wh}) 
most remarkably the case of stationary varifolds or area minimizing currents
in a Riemannian manifold.
For a more complete account on the theory of varifolds and almost minimizing currents
we refer to \cite{Allard}, \cite{Bomb} and the lecture notes \cite{Sim83}.

\subsection{Tubular neighborhood estimate}\label{s:vrfld tub neg}
In what follows we consider integer rectifiable varifolds $\sV=(\Gamma,f)$,
where $\Gamma$ is an $m$-dimensional rectifiable set in the bounded open subset
$\Omega\subset\R^n$, and $f:\Gamma\to\N\setminus\{0\}$ is locally $\cH^m$-integrable.
We assume that $\sV$ has bounded generalized mean curvature, i.e.~there exists a vector
field $H_\sV : \Omega \to \R^n$ such that $\|H_\sV\|_{L^\infty(\Omega,\R^n)} \leq H_0$
for some $H_0 >0$ and
\[
\int_\Gamma\mathrm{div}_{T_y\Gamma} X\,d\mu_\sV = 
-\int X\cdot H_\sV\,d\mu_\sV
\qquad\forall\; X\in C^1_c(\Omega,\R^n)
\]
where $\mu_\sV := f \, \cH^m \res \Gamma$.
It is then well-known (cp., for example, \cite[Theorem~17.6]{Sim83}) that the quantity
\[
\Theta_\sV(x,\rho):=e^{H_0\,\rho}\,\frac{\mu_\sV(B_\rho(x))}{\omega_m\rho^m}
\]
is monotone and the following inequality holds for all $0<\sigma<\rho<\mathrm{dist}(x,\partial\Omega)$
\begin{equation}\label{e:monotineq}
\Theta_\sV(x,\rho)-\Theta_\sV(x,\sigma)\geq
\int_{B_\rho\setminus B_\sigma(x)}\frac{|(y-x)^\perp|^2}{|y-x|^{m+2}}d\mu_\sV(y)
\end{equation}
where $(y-x)^\perp$ is the orthogonal projection of $y-x$ on the orthogonal complement $(T_y\Gamma)^\perp$.
In particular the family $(\Theta(\cdot,s))_{s\in [0,r_0]}$ (with the obvious
extended notation $\Theta(\cdot,0^+):= \lim_{r\downarrow 0}\Theta(\cdot, r)$)
satisfies assumption (a) in \S~\ref{pa:1.1} for every fixed $r_0>0$ with
\begin{equation}\label{e:lambdavar}
\Lambda_0(r_0):=e^{H_0\,\mathrm{diam}(\Omega)}\,\frac{\mu_{\sV}(\Omega)}{\omega_m r_0^{m}}.
\end{equation}

\medskip

In order to introduce the control functions $\d_k$ we recall next the definition of \textit{cone}.

\begin{definition}\label{d.cone}
An integer rectifiable $m$-varifold $\sC = (R, g)$ in $\R^n$ is a cone if
the $m$-dimensional rectifiable set $R$ is invariant under dilations 
i.e.
\[
\lambda\,y \in R \qquad \forall \;y \in R, \;\forall\;\lambda>0 
\]
and $g$ is $0$-homogeneous, i.e.
\[
g(\lambda\,y) = g(y) \qquad \forall\; y \in R,\; \forall \lambda>0. 
\]
An integer rectifiable $m$-varifold $\sC = (R, g)$ in $B_\rho$, $\rho>0$ is a cone if
it is the restriction to $B_\rho$ of a cone in $\R^n$.

The spine of a cone $\sC = (R, g)$ in $\R^n$ is the biggest subspace $V \subset \R^n$
such that $R = R' \times V$ up to $\cH^m$-null sets. 

The class of cones whose spine is at least $k$-dimensional is
denoted by $\mathcal{C}_k$ and its elements are called $k$-conical.
\end{definition}

If $\d_\ast$ is a distance inducing the weak $\ast$ topology of varifolds
with bounded mass in $B_1$
(cp., for instance, \cite[Theorem~3.16]{Rud91} for the general case of dual spaces),\
the control function $\d_k$ is then defined as
\begin{equation}\label{e:dkvar}
\d_k(x,s):=\inf\big\{\d_\ast\big(\sV_{x,s},\sC\big):\,\sC \in \mathcal{C}_k, \; \|H_{\sC}\|_{L^\infty(\Omega,\R^n)}\leq H_0\big\}
\end{equation}
where $\sV_{x,s}:=(\eta_{x,s}(\Gamma),f\circ \eta_{x,s}^{-1})$ with
$\eta_{x,s}(y):=\sfrac{(y-x)}{s}$.

By very definition, then (b) in \S~\ref{pa:1.1} is satisfied.
We are now ready to check that the conditions in the Structural Hypotheses are satisfied. 
As usual, we write the corresponding statements for fixed $r_0$ and $\Lambda_0:=\Lambda_0(r_0)$, for simplicity.

\begin{lem}\label{l:qrvar}
For every $\eps_1>0$ there exist $0<\lambda_1(\eps_1),\,\eta_1(\eps_1) <\sfrac{1}{4}$
such that for all $(x,\rho) \in U$, with $x\in\Omega^{r_0}$ and $\rho<r_0$,
\[
\Theta_\sV(x,\rho)-\Theta_\sV(x,\lambda_1\,\rho)\leq\eta_1
\quad\Longrightarrow\quad
\d_0(x,\rho) \leq \eps_1.
\]
\end{lem}

\begin{pf}
Assume by contradiction that for some $\eps_1>0$ there exists
$(x_j,\rho_j)\in U$, with $x_j\in\Omega^{r_0}$ and $\rho_j<r_0$, such that
\begin{equation}\label{e:contravar}
\Theta_{\sV}(x_j,\rho_j)-\Theta_{\sV}(x_j,j^{-1}\,\rho_j)\leq j^{-1}
\quad\text{and}\quad\d_0(x_j,\rho_j) \geq \eps_1.
\end{equation}
We consider the sequence $(\sV_j)_{j \in \N}$ with $\sV_j:=\sV_{x_j,\rho_j}$,
and note that for all positive
$t\,>0$ there is an index $\bar{j}$ such that $t\,\rho_j<r_0$ if $j\geq\bar{j}$,
so that
\[
\mu_{\sV_j}(B_t\,)\leq\omega_m\,t^m\,\Theta_{\sV_j}(x_j,t\,\rho_j)\leq 
 \omega_mt\,^m\Lambda_0\quad\forall\; j\geq\bar{j}.
\]
Therefore, up to the extraction of subsequences and a diagonal argument,
Allard's rectifiability criterion (cp., for instance, \cite[Theorem~42.7, Remark 42.8]{Sim83}) yields a limiting $m$-dimensional 
integer varifold $\sV_j \to \sC=(R,g)$ with the bound
$\|H_\sC\|_{L^\infty(\Omega} \leq H_0$.
Since $\Theta_\sV(x_j,s\,\rho_j)=\Theta_{\sV_j}(0,s)
\to\Theta_\sC(0,s)$ except at most for countable values of $s$,
by monotonicity and \eqref{e:contravar} for all $j^{-1}<r<s<1$
we have $\Theta_\sC(0,s) = \Theta_\sC(0,0^+)$ for every $s\geq 0$.
The monotonicity formula \eqref{e:monotineq} applied to $\sC$ implies that $\sC$
is actually a cone, thus contradicting $\d_0(x_j,\rho_j)\leq \eps_1$.
\end{pf}

\begin{lem}\label{l:csvar}
For every $\eps_2,\tau\in (0,1)$, there exists $0<\eta_2(\eps_2,\tau)<\eps_2$
such that, for every $(x,5s)\in U$, with $x\in\Omega^{r_0}$ and $5s<r_0$, if for some $k \in \{0, \ldots, m-1\}$
\[
\d_k(x,4s) \leq \eta_2\quad\text{ and }\quad
\d_{k+1}(x,4s) \geq \eps_2,
\]
then there exists a $k$-dimensional affine space $x+V$ such that
\[
\d_0(y,4s)> \eta_2 \quad \forall\; y \in B_{s}(x)\setminus
\cT_{\tau s}(x+V).
\]
\end{lem}

\begin{pf}
The proof is by contradiction. Assume that there exist
$0<\eps_2,\,\tau<1$, $k \in \{0, \ldots, m-1\}$ and a sequence of points
$(x_j,5s_j)\in U$, with $x_j\in\Omega^{r_0}$ and $5s_j<r_0$, for $2j \geq \eps_2^{-1}$ such that
\begin{equation}\label{e:cono1}
\d_{k}(x_j,4s_j)\leq j^{-1}\quad\text{ and }\quad \d_{k+1}(x_j,4s_j) \geq \eps_2,
\end{equation}
and such that the conclusion of the lemma fails, in particular, for $V_j$ given by the spine of $\sC_j$ with
\begin{equation}\label{e:vicinanza cono invariante}
\d_\ast\big(\sV_{x_j,4s_j},\sC_j\big)\leq 2j^{-1}
\end{equation}
(note that by $2j \geq \eps_2^{-1}$ necessarily $\dim(V_j) = k$).
Without loss of generality (up to a rotation) we can assume that $V_j = V$ a given
vector subspace for every $j$.
This means that there exist $y_{j} \in B_{s_j}(x_j)\setminus\cT_{\tau s_j}(x_j+V)$ such that
\begin{equation}\label{e:cono2}
\d_0(y_{j},4s_j)\leq j^{-1}. 
\end{equation}
Using the compactness for varifolds with bounded generalized mean curvature,
(up to passing to subsequences) we can assume that 
\begin{enumerate}
\item $s_j \to s_\infty \in [0, \sfrac{r_0}{5}]$;
\item $\sC_j \to \sC_\infty$ in the sense of varifolds, $\sC_\infty$ a cone with
$\|H_{\sC_\infty}\|_{L^\infty(\Omega,\R^n)} \leq H_0$;
\item $\sfrac{(y_j-x_j)}{s_j}\to z\in \overline{B}_{1}\setminus\cT_{\tau}(V)$;
\item $\sV_{x_j,s_j} \to \sW_\infty$ and $\sV_{y_j,s_j} \to \sZ_\infty$ in the
ball $B_4$ in the sense of varifolds, 
where $\sW_\infty$ and $\sZ_\infty$ are cones thanks to
\eqref{e:cono1} and \eqref{e:cono2}, respectively.
\end{enumerate}
Note that by \eqref{e:vicinanza cono invariante} it follows that $\sC_j \to \sW_\infty$
and therefore $\sW_\infty \in \mathcal{C}_k$ because all the $\sC_j$ are invariant
under translations in the directions of $V$.
Moreover, arguing as above it also follows from $\d_{k+1}(x_j,4s_j) \geq \eps_2$
that the spine of $\sW_\infty$ is exactly $V$.

Note that $\eta_{\sfrac{(y_j - x_j)}{s_j},1}$ corresponds to the translation of vector
$\sfrac{(y_j - x_j)}{s_j}$. By the equality of $(\eta_{\sfrac{(y_j - x_j)}{s_j},1})_\sharp \sV_{x_j, s_j}$
and $\sV_{y_j,s_j}$ in $B_3$,
we deduce that $\sZ_\infty=(\eta_{\sfrac{(y_j - x_j)}{s_j},1})_\sharp \sW_\infty$
as varifolds in $B_3$, i.e.~$\sW_\infty$ is a cone around $z$ too.
We claim that this implies that $\sW_\infty$ is invariant along the directions of $\mathrm{Span}\{z, V\}$, 
thus contradiction the fact that the spine of $\sW_\infty$ equals $V$.
To prove the claim, let $\sW_\infty = (R_\infty, g)$ with $R_\infty$ cone around the origin and $z$. It 
suffices to show that $y+z\in R_\infty$ for all $y\in R_\infty$.
Indeed $\sfrac{(z+y)}{2} = z + \sfrac{y-z}{2}\in R_\infty$ being $R_\infty$ a cone with respect to $z$; and then $y+z \in R_\infty$ being $R_\infty$ a cone with respect to $0$.
\end{pf}

In particular we deduce that Theorem~\ref{t:tub-neig} and
Theorem~\ref{t:nostra stratiFICA}
hold in the case of varifolds with generalized mean curvature in $L^\infty$.

\subsection{Almost minimizer in codimension one}
It is well-known by the classical examples by Federer \cite{Fed} that no
Allard's type $\eps$-regularity results can hold
for higher codimension generalized submanifolds
without any extra-hypotheses on the densities.
Vice versa for generalized hypersurfaces one can strengthen the results of the previous
subsection giving estimates on the Minkowski dimension of the singular set.
The arguments in this part resemble very closely those in \cite{ChHaNa2},
therefore we keep them to the minimum.

In what follows we consider sets of finite perimeter, i.e.~borel subsets
$E \in \Omega$ such that the distributional derivative of corresponding
characteristic function has bounded variation: $D\chi_E \in BV_\Omega$.
Following \cite{Bomb,Tam84}, a set of finite perimeter is almost minimizing in $\Omega$ if for all $A\subset\subset\Omega$ open there exist $T\in\big(0,\dist(A,\partial\Omega)\big)$ and $\alpha:(0,T)\to[0,+\infty)$ non-decreasing and infinitesimal in $0$ such that whenever 
$E\triangle F\subset\subset B_r(x)\subset A$
\begin{equation}
\mathrm{Per}(E,B_r(x))\leq\mathrm{Per}(F,B_r(x))+\alpha(r)\,r^{n-1}\qquad\forall\; r\in(0,T)\label{e:am}
\end{equation}
and
\begin{equation}
(0,T)\ni t\mapsto \frac{\alpha(t)}t\,\text{ is non-increasing, and }\,
 \int_0^T\frac{\alpha^{\sfrac{1}{2}}(t)}tdt<\infty.\label{e:am2}
\end{equation}

Examples of almost minimizing sets not only include minimal boundaries on Riemannian manifolds, 
but also boundaries with generalized mean curvature in $L^\infty$, minimal boundaries with volume constraint, and minimal 
boundaries with obstacles (cp. \cite[\S~1.14]{Tam84}).

We use here again the control functions introduced in Section~\ref{pa.mc}
in terms of flat distance: given a set of finite perimeter $E$, we denote by $\de E$
its boundary (in the sense of currents) and set
\[
\d_{k}(x,s) := \inf\big\{\F\big((\de E_{x,s} - C) \res B_1\big) \;:\; C \; \text{ $k$-conical \& area minimizing} \big\}
\]
where the dimension of the cones $C$ is always $n-1$,
and $E_{x,s}$ is the push-forward of $E$ via the rescaling map $\eta_{x,s}$.
In particular $d_{n-1}$ denotes the distance of the rescaled boundary $\de E_{x,s}$
rescaling of the from flat $(n-1)$-dimensional vector spaces.

The main $\eps$-regularity result for almost minimizing sets can be stated as follows
(cp.~\cite[Theorem~1.9]{Tam84}, \cite[Lemma 17]{Bomb} and \cite[Theorem~B.2]{Sim83}).

\begin{thm}\label{t:Tam}
Suppose that $E$ is a perimeter almost minimizer in $\Omega$ satisfying
\eqref{e:am} and \eqref{e:am2} for a given function $\alpha$.
Then, there exists $\varepsilon>0$ and
$\omega:[0,+\infty) \to [0,+\infty)$ continuous, non-decreasing and satisfying $\omega(0)=0$ 
with the following property: if
\[
\rho+\d_{n-1}(x,\rho)+\int_0^\rho \frac{\alpha^{\sfrac{1}{2}}(t)}tdt\leq\varepsilon,
\]
then $\de E \cap B_{\sfrac{\rho}{2}}(x)$ is the graph of a $C^{1}$ function
$f$ satisfying
\begin{equation}\label{e:uniform cont}
|\nabla f(x) - \nabla f(y)| \leq \omega(|x-y|).
\end{equation}
Moreover, there are no singular area minimizing cones with dimension of the
singular set bigger than $n-8$, i.e.~equivalently
\begin{equation}\label{e:no sing cones}
\d_{n-7} = \d_{n-6} = \ldots = \d_{n-1}.
\end{equation}
\end{thm}

\begin{rmk}
The smallness condition $\d_{n-1} \leq \eps$, together with the
almost minimizing property, implies 
the more familiar smallness condition on the \textit{Excess}, i.e.
\[
\textup{Exc}(E,B_r(x)):=r^{1-n}\,\|D\chi_E\|(B_r(x))-r^{1-n}
\left|D\chi_E (B_r(x))\right| \leq \varepsilon'
\]
for some $\varepsilon' = \varepsilon'(\varepsilon)>0$ infinitesimal
as $\varepsilon$ goes to $0$
because of the continuity of the mass for converging uniform almost
minimizing currents.
Therefore \eqref{e:uniform cont} readily follows from \cite[Theorem~1.9]{Tam84}.
\end{rmk}

By a simple use of Theorem~\ref{t:Tam} we can the prove the following.

\begin{corollary}\label{c:eps reg}
Under the hypotheses of Theorem~\ref{t:Tam} there exist constants $\delta_0 = \delta_0(\Lambda_0,n, \alpha)>0$ 
and $\rho_0=\rho_0(\Lambda_0, n, \alpha) >0$ such that
\[
\mathcal{S}^{n-8}_{r_0,\delta_0}= \mathcal{S}^{n-8}_{r_0} = 
\mathcal{S}^{n-7}_{r_0} = \ldots = \mathcal{S}^{n-2}_{r_0}
\qquad \forall \; r_0 \in (0,\rho_0].
\]
\end{corollary}

\begin{pf}
Set $\delta_0 = \sfrac{\eps}{2}$ and let $\rho_0$ be sufficiently small to have
\[
\rho_0
+\int_0^{\rho_0} \frac{\alpha^{\sfrac{1}{2}}(t)}tdt\leq\sfrac{\varepsilon}{2}.
\]
If $x \not\in \mathcal{S}^{n-2}_{r_0,\delta_0}$, $r_0\in(0,\rho_0]$, 
then there exists $0<z_0 \leq r_0$ such that $\d_{n-1}(x,z_0) < \delta_0$.
In particular, by the choices of $\delta_0$ and of $\rho_0$ the assumptions of Theorem~\ref{t:Tam}
are satisfied at $s_0$. Therefore, it turns out that $x$ is a regular point of $\de E$ and that 
$B_{\sfrac{z_0}2}(x) \cap \de E$ can be written as a graph of a function $f$ satisfying \eqref{e:uniform cont}. 
In particular, $\lim_{s\downarrow 0} \d_{n-1}(x,s) =0$.
Therefore, given any $\delta'<\delta_0$, we have that $x \not\in \mathcal{S}^{n-2}_{r_0,\delta'}$, 
thus implying that $\mathcal{S}^{n-2}_{r_0}=\mathcal{S}^{n-2}_{r_0,\delta_0}$. 
By taking into account \eqref{e:no sing cones} we conclude the corollary straightforwardly. 
\end{pf}

In particular, Theorem~\ref{t:minchioschi} holds and we deduce the following
refinement of the Hausdorff dimension estimate of the singular set.

\begin{thm}\label{t:am minchioschi}
Let $E \subset \Omega$ be a almost minimizing set of finite perimeter in a bounded
open set $\Omega\subset \R^n$ according to \eqref{e:am} and \eqref{e:am2}.
Then there exists a closed subset $\Sigma \subset \de E \cap \Omega$ such that 
$\de E \cap \Omega \setminus \Sigma$ is a $C^1$ regular $(n-1)$-dimensional submanifold 
of $\R^n$ and $\dim_\cM(\Sigma) \leq n-8$.
\end{thm} 

\begin{pf}
Let $\Omega' \subset \subset \Omega$ be compactly supported and set $r_0 := \dist(\Omega', \de \Omega)$.
By the regularity Theorem~\ref{t:Tam}, a point $x \in \Omega$ is regular
if and only if there exists $r>0$ sufficiently small such that $\d_{n-1}(x,r) \leq \sfrac{\eps}{2}$. 
In particular, the set of singular points $\Sigma$ coincides with
$\mathcal{S}^{n-2}_{r_0,\sfrac{\eps}{2}}$ and the conclusion follows combining
Theorem~\ref{t:minchioschi} with Corollary~\ref{c:eps reg}.
\end{pf}

In addition, we can also derive a higher integrability estimate for almost minimizers with bounded
generalized mean curvature.
Given a set of finite perimeter $E\subset \Omega$, one can associate to $\de E$
a varifold in a canonical way (cp.~\cite{Sim83}). One can then talk about
sets of finite perimeter with bounded generalized mean curvature.
Important examples of such an instance are:
\begin{enumerate}
\item the minimizers of the area functional in a Riemannian manifold;
\item the minimizers of the prescribed curvature functional in $\Omega \subset \R^n$
\[
\mathcal{F}(E) := \|D\chi_E\|(\Omega) + \int_{\Omega \cap E} H
\]
with $H \in L^\infty(\Omega)$;
\item minimizers of the area functional with volume constraint;
\item more general $\Lambda$-minimizers for some $\Lambda>0$, i.e.~sets $E$ such that
\[
\|D\chi_E\|(\Omega) \leq \|D\chi_F\|(\Omega) + \Lambda\,|E \setminus F|
\quad \forall \; F \subset\Omega.
\]
\end{enumerate}
Given a point $x \in \de E$ such that $B_r(x) \cap \de E$ is the graph of a $C^1$
function $f$, if the generalized mean curvature $H$ of $\de E$ is bounded then we can 
also talk about generalized second fundamental form $A$ in $B_{\sfrac{r}{2}}(x)$, because in a suitable chosen system of coordinates $f$ solves in a weak sense the 
prescribed mean curvature equation
\begin{equation}\label{e:equazione minchiona}
\textup{div} \left( \frac{\nabla f}{\sqrt{1+ |\nabla f|^2}}\right) = H \in L^\infty.
\end{equation}
Note that, since in this case $f$ satisfies \eqref{e:uniform cont}, we can 
choose a suitable system of coordinates and use
the $L^p$ theory for uniformly elliptic equations to deduce that
actually $A \in L^p(B_{\sfrac{r}{4}}(x), \cH^{n-1}\res \de E)$ for every $p< +\infty$
with uniform estimate
\begin{equation}\label{e:Lp ellittico}
\int_{B_{\frac{r}{4}}(x) \cap \de E} |A|^p \, \cH^{n-1} \leq C\, r^{n-p-1}
\end{equation}
for some dimensional constant $C>0$.
For convenience we set $A \equiv +\infty$ on the singular set $\Sigma \subset \de E$.

\begin{thm}\label{t:higher int}
Let $E \subset \Omega$ be as in Theorem~\ref{t:am minchioschi} and
assume moreover that the varifold induced by $\de E$ has bounded generalized mean curvature.
Then, for every $p<7$ there exists a constant $C>0$ such that
\begin{equation}\label{e:higher int}
\int_{\de E \cap \Omega} |A|^p \,\d\cH^{n-1} \leq C.
\end{equation}
\end{thm}

\begin{pf}
Let $\rho_0>0$ be the constant in Corollary~\ref{c:eps reg} and $\varepsilon$ that
of Theorem~\ref{t:Tam}.
Then $\Sigma = \mathcal{S}^{n-8}_{\rho_0, \sfrac{\varepsilon}{2}}$.
In then follows that for a fixed $\bar k > \log_2(\sfrac{\rho_0}{10})$
\[
\big(\supp(\de E) \setminus \Sigma \big) \cap \Omega 
= \bigcup_{k \geq \bar k} \mathcal{S}^{n-8}_{2^{-k}, \rho_0, \sfrac{\eps}{2}} \setminus
\mathcal{S}^{n-8}_{2^{-k-1}, \rho_0, \sfrac{\eps}{2}}.
\]
Applying Theorem~\ref{t:tub-neig} we infer that for every $\eta>0$ there exists
$C>0$ such that
\begin{equation}\label{e:stima volume}
\big\vert \mathcal{T}_{2^{-k}}(\mathcal{S}^{n-8}_{2^{-k}, \rho_0, \sfrac{\eps}{2}}) \big\vert \leq C\, 2^{-k (8 - \eta)}.
\end{equation}
By Lemma~\ref{l.ricoprimento banana} there exists a cover of
$\mathcal{T}_{\sfrac{2^{-k-2}}{5}}(\mathcal{S}^{n-8}_{2^{-k}, \rho_0, \sfrac{\eps}{2}} \setminus \mathcal{S}^{n-8}_{2^{-k-1}, \rho_0, \sfrac{\eps}{2}})$
by balls $\{B_{2^{-k-3}}(x^k_i)\}_{i \in I_k}$ with $x^k_i \in \mathcal{S}^{n-8}_{2^{-k}, \rho_0, \sfrac{\eps}{2}} \setminus \mathcal{S}^{n-8}_{2^{-k-1}, 
\rho_0, \sfrac{\eps}{2}}$ whose cardinality is estimated by \eqref{e.stima banana} as
\begin{equation}\label{e:stima cardinalita'}
\cH^0(I_k) \leq C\, 2^{-k (8 - \eta - n)}
\end{equation}
where $C>0$ is a dimensional constant.

\medskip

We start estimating the integral in \eqref{e:higher int} as follows:
\begin{align}
\int_{\de E \cap \Omega} |A|^p \,\d\cH^{n-1} & = \sum_{k \geq \bar k} \int_{\mathcal{S}^{n-8}_{2^{-k}, \rho_0, \sfrac{\eps}{2}} \setminus
\mathcal{S}^{n-8}_{2^{-k-1}, \rho_0, \sfrac{\eps}{2}}} |A|^p \,\d\cH^{n-1}\notag\\
& \leq \sum_{k\geq \bar k} \sum_{i \in I_k} \int_{\de E \cap B_{2^{-k-3}}(x^k_i)} |A|^p\,\d\cH^{n-1}\notag
\end{align}
Since $x^k_i \in \mathcal{S}^{n-8}_{2^{-k}, \rho_0, \sfrac{\eps}{2}} \setminus \mathcal{S}^{n-8}_{2^{-k}, \rho_0, \sfrac{\eps}{2}}$ it follows that there exists
$r_i^k \in [2^{-k-1}, 2^{-k})$ such that $d_n(x^k_i,r^k_i) < \sfrac{\eps}{2}$.
In particular by Theorem~\ref{t:Tam} $\de E \cap B_{2^{-k-2}}(x^k_i)$ is a graph
of a $C^1$ function satisfying \eqref{e:uniform cont}.
From \eqref{e:Lp ellittico} we conclude that 
\begin{align}
\int_{\de E \cap \Omega} |A|^p \,\d\cH^{n-1} & \leq C \sum_{k \geq \bar k} \cH^0(I_k)\, 2^{-k(n-p-1)} \leq C \sum_{k \geq \bar k} 2^{-k(7 - \eta -p)} < C\notag
\end{align}
as soon as $\eta< 7 - p$.
\end{pf}

%
%

%

\bibliography{Q-strat.bib}

\end{document}